\documentclass{article}

\usepackage{amsmath}
\usepackage{amssymb}
\usepackage{latexsym}
\usepackage{amsthm}
\usepackage{Young}

\theoremstyle{plain} \newtheorem{cor}{Corollary}[section]
\theoremstyle{plain} \newtheorem{lem}[cor]{Lemma} 
\theoremstyle{plain}  
\theoremstyle{plain} \newtheorem{thm}[cor]{Theorem}
\theoremstyle{plain} \newtheorem*{thma}{Theorem}
\theoremstyle{plain}  
\theoremstyle{plain} \newtheorem*{defa}{Definition}
\theoremstyle{plain} \newtheorem*{lema}{Lemma} 
\theoremstyle{plain}  
\theoremstyle{definition} \newtheorem{consideration}{Consideration}
\theoremstyle{definition} \newtheorem{rmk}[cor]{Remark}

\newcommand{\dom}{\textup{dom}}
\newcommand{\ran}{\textup{ran}}
\newcommand{\rk}{\textup{rk}}
\newcommand{\ld}{\left\lfloor}
\newcommand{\rd}{\right\rfloor}
\newcommand{\Z}{\mathbb Z}
\newcommand{\ZnZ}{\Z / n\Z}
\newcommand{\C}{\mathbb C}
\newcommand{\CS}{\mathbb C S}
\newcommand{\CRn}{\mathbb C R_n}
\newcommand{\Y}{\mathcal Y}
\newcommand{\s}{\ld s \rd}
\newcommand{\onetok}{\{1,\ldots,k\}}
\newcommand{\zeroton}{\{0,\ldots,n\}}
\newcommand{\oneton}{\{1,\ldots,n\}}
\newcommand{\semigroupcplx}{{\mathcal T}_{\Y}^s}
\newcommand{\groupoidcplx}{{\mathcal T}_{\Y}^{\s}}
\newcommand{\minup}{\textup{min}}
\newcommand{\hyph}{-\penalty0 \hskip0pt\relax }
\newcommand{\LP}{\textup{(}}
\newcommand{\RP}{\textup{)}}

\title
{Fast Fourier Transforms for the Rook Monoid}

\author{Martin Malandro and Dan Rockmore\footnote{supported by AFOSR under grant FA9550-06-1-0027.}}

\begin{document}






\maketitle						


\begin{abstract}
We define the notion of the Fourier transform for the rook monoid (also called the symmetric inverse semigroup) and provide two efficient divide-and-conquer algorithms (fast Fourier transforms, or FFTs) for computing it. This paper marks the first extension of group FFTs to non-group semigroups.
\end{abstract}

\section{Introduction}

The rook monoid $R_n$, also called the symmetric inverse semigroup, is the set of all partial permutations of $\{1,\ldots,n\}$ under function composition, i.e., it is the semigroup of all bijections between all subsets of $\{1,\ldots,n\}$ under function composition, with the usual rule for composing partial functions. That is, $g\circ f$ is defined for precisely the elements $x$ for which $x \in \dom (f)$ and $f(x) \in \dom (g)$. It is called the rook monoid because it is isomorphic to the semigroup of all $n \times n$ matrices with the property that at most one entry in each row is 1 and at most one entry in each column is 1 (the rest being 0) under multiplication. Such matrices (called {\em rook matrices}) correspond to the set of all possible placements of non-attacking rooks on an $n \times n$ chessboard. For example, consider the element $\sigma \in R_4$ defined by
\[\sigma(2)=1, \quad \sigma(4)=4.             \]
Then, viewed as a partial permutation, $\sigma$ is
\begin{displaymath}
\sigma = \left( \begin{array}{cccc}
1&2&3&4\\
-&1&-&4 \end{array} \right)
\end{displaymath}
where the dash indicates that the above entry is not mapped to anything. As a rook matrix, we have
\begin{displaymath}
\sigma = \left[
\begin {array}{cccc}
0&1&0&0\\
0&0&0&0\\
0&0&0&0\\
0&0&0&1
\end {array}
\right]. \end{displaymath}

In this paper we define the notion of the {\em Fourier transform} for a $\C$-valued function on the rook monoid, and we address the problem of computing it efficiently. We present two distinct fast Fourier transform (FFT) algorithms, both of which make use of tools generalized from group FFT theory \cite{Dan2004}. We assume no specialized knowledge. This paper is meant to be readable by group FFT theorists with an interest in semigroups as well as semigroup theorists with an interest in FFTs.

The history of the FFT is an interesting one, beginning, like many subjects in mathematics, with Gauss, who developed the now-classical FFT to efficiently interpolate the orbits of asteroids. This FFT was rediscovered by Cooley and Tukey \cite{CooleyTukey} in 1965 to analyze seismic data. 
A different FFT was discovered in 1937 by Yates \cite{Yates}, for the analysis of data on factorial designs. The theory of FFTs for Abelian groups reconciles these two FFTs. The first is an FFT on $\mathbb Z/n\mathbb Z$, and the second is an FFT on $(\mathbb Z/2\mathbb Z)^k$. From this generalization also grew a theory of FFTs for non-Abelian groups. For example, we now have FFTs on supersolvable groups, FFTs on the symmetric group, FFTs on compact Lie groups, and FFTs on finite groups of Lie type. For a survey of these results, see \cite{Dan2004}.

We have both theoretical and practical motivations for creating FFT algorithms for the rook monoid. From a theoretical perspective, these FFTs are significant because they are the first examples of FFTs on a non-group semigroup. The rook monoid plays the same role for finite inverse semigroups that the symmetric group $S_n$ does for finite groups in terms of Cayley's theorem: every finite inverse semigroup is isomorphic to a sub-semigroup of some rook monoid (see, e.g. \cite{Lawson}, p. 36-37), so this makes the rook monoid a natural place to begin extending the theory of FFTs to semigroups. Also, just as groups capture symmetry, inverse semigroups capture partial symmetry. See [16] for more on this idea. 


From a practical perspective, partially ranked data occurs naturally, and the rook monoid provides a nice computational framework for analyzing such data. For example, consider voting data in which some voters do not fully rank all candidates. In practice, this may well be the case. For example, in the 1980 American Psychological Association (APA) presidential election, only a third of the voters cast fully ranked ballots \cite{Persi}. Say $f(\sigma)$ is the number of voters who submitted a ballot of type $\sigma$. For the $\sigma$ above, $f(\sigma)$ would represent the number of voters who ranked candidate 2 in the first position, candidate 4 in the fourth position, and didn't bother ranking candidates 1 and 3. The values $\{f(\sigma)\}$ make up the dataset. There have been previous attempts at spectral analysis for such datasets. For example, in \cite{Persi}, Diaconis used spectral analysis on symmetric groups to analyze the APA election mentioned above. We propose viewing $f$ as an element of the semigroup algebra $\mathbb C R_n$, so that spectral analysis of $f$ involves decomposing $\mathbb C R_n$ into its isotypic components by means of an FFT. Things are more complicated than in the group case because $\mathbb C R_n$ has two natural bases. This means that there are two natural ways to view functions on $R_n$ as elements of $\mathbb C R_n$, and therefore two different notions for convolution of such functions. We will address both in this paper. 

Partially ranked voting data is also a type of ``not missing at random" (NMAR) data. Missing data in a sample is said to be NMAR if it is believed to be missing, at least in part, because of its unobserved value. In the voting example, any ballot that is not fully ranked is considered NMAR, since the voter intentionally left something on the ballot blank. Consumer survey data frequently contains NMAR data. NMAR data contrasts with randomly missing data, which, in the voting example, can be caused by accidentally losing some of the ballots. To date, there doesn't seem to be a standard method for analyzing NMAR data. However, the concept of spectral analysis for the rook monoid immediately applies to any NMAR dataset that can be considered a collection of partial rankings. Since spectral analysis is a model-independent approach, spectral analysis for the rook monoid may go a long way towards creating a standard method for dealing with NMAR data.


This paper is organized as follows. Sections \ref{Defs} and \ref{RnFacts} contain basic definitions and background material about semigroups and the rook monoid. Section \ref{SemigpReps} contains a brief introduction to semigroup representation theory. Sections \ref{AdaptedReps} and \ref{Schurs} are devoted, respectively, to two important tools in FFT theory: chain-adapted matrix representations and Schur's Lemma, adapted to semigroups.

Section \ref{BasisConsid} is devoted to basis considerations for semigroup algebras. Section \ref{NaturalBases} defines both natural bases for $\C R_n$, and Section \ref{FourierBases} defines the notion of a Fourier basis. Section \ref{FourierTransformCS} defines the Fourier transform.

Section \ref{SteinDecomp} shows how the semigroup algebra $\C R_n$ decomposes into a product of matrix algebras over group algebras. This has far-reaching implications, some of which are explored in Section \ref{SteinDecompConseq}. Section \ref{SteinReps} describes explicit matrix representations for $R_n$, Section \ref{ExplicitFourierBasis} contains an explicit Fourier basis description for $\C R_n$, and Section \ref{FourierInversion} contains the Fourier inversion theorem for $\C R_n$.

In Section \ref{FunctionAssoc}, we explain some of the issues that arise when choosing how to associate functions on $R_n$ to elements of $\C R_n$. As mentioned previously, this is only an issue because $\C R_n$ has two natural bases.

Section \ref{FFTConsiderations} contains FFT-specific considerations. Computational complexity for FFTs is defined in Section \ref{CompComplexity}. Our algorithms for the FFT on $R_n$ use ideas from group FFT theory, which we cover in Sections \ref{FFTsOnGroups} and \ref{SnFFT}.

Section \ref{RnFFT1} contains our better FFT algorithm for $R_n$. It relies heavily on the decomposition described in Section \ref{SteinDecomp} and the results in Sections \ref{SteinReps} and \ref{SnFFT}.

Section \ref{RnFFT2} contains our other FFT algorithm for $R_n$. The ideas involved are mainly from Sections \ref{AdaptedReps}, \ref{Schurs}, and \ref{FFTsOnGroups}. While it is not as efficient as the algorithm given in Section \ref{RnFFT1}, it is constructed in an entirely different way.  We present this algorithm because the ideas involved may be useful for constructing FFTs on other semigroups.


\section{Preliminaries}
\label{Preliminaries}

\subsection{Definitions for general semigroups}
\label{Defs}

\begin{defa}[semigroup]
A {\em semigroup} is a nonempty set $S$ together with an associative binary operation \textup{(}which we write multiplicatively\textup{)}.
\end{defa}

\begin{defa}[inverse semigroup]
An {\em inverse semigroup} is a semigroup $S$ such that, for each $x \in S$, there is a {\em unique} $y \in S$ such that
\[xyx = x \textup{ and } yxy = y. \] In this case, we write $y=x^{-1}$.
\end{defa}

We remark that the condition that $y$ be unique is necessary for this definition. An element $x \in S$ is said to be {\em regular} or {\em Von-Neumann regular} if there is at least one $y \in S$ satisfying $xyx = x$ and $yxy = y$, and $S$ is said to be {\em regular} if every element of $S$ is regular. Consider the full transformation semigroup $X$ on the set $\{1, 2, \ldots, n\}$; that is, all maps from $\{1, 2, \ldots, n\}$ to itself under composition. It is easy to see that $X$ is regular, and that (for $n \geq 2$) there exist elements $x \in X$ for which there are {\em multiple} elements $y \in X$ satisfying $xyx = x$ and $yxy = y$. $X$ is therefore not inverse.  An equivalent characterization of inverse semigroups (see, e.g. \cite{Lawson}) is as follows.

\begin{defa}[inverse semigroup] An {\em inverse semigroup} is a semigroup $S$ which is regular and for which all idempotents of $S$ commute.\end{defa}

\begin{defa}[semigroup algebra] Let $S$ be a finite semigroup. The {\em semigroup algebra} $\mathbb{C}S$ is the formal $\mathbb C$-span of the symbols $\{s\}_{s\in S}$. Multiplication in $\mathbb C S$, denoted by $\ast$, is given by convolution \LP i.e., the linear extension of the semigroup operation via the distributive law\RP : Suppose $f, g \in \mathbb C S$, with
\[f = \sum_{r\in S}f(r) r, \quad g=\sum_{t\in S}g(t) t.\]
Then
\begin{equation}
f\ast g = \sum_{r\in S}f(r) r \sum_{t\in S}g(t) t = \sum_{s \in S} \sum_{r,t \in S: rt=s} f(r) g(t) s.
\label{CSConvolutionSemigpBasis}
\end{equation}
\end{defa}

If $S$ is a group, then convolution may be written in the familiar way:

\[f\ast g = \sum_{s \in S}\sum_{r \in S} f(r)g(r^{-1}s)s.\]


Note that convolution here means convolution {\em in the semigroup algebra}. We hold off on defining a notion of convolution for {\em functions} on $S$ until Section \ref{FunctionAssoc}. We will use the semigroup algebra convolution to define function convolution once we have defined how to associate functions to elements of the semigroup algebra. It turns out that, unlike in the case of groups, there is in general more than one natural way to associate functions on $S$ to elements of $\mathbb C S$ because, in general, $\mathbb C S$ has more than one natural basis.

\subsection{Facts about $R_n$}
\label{RnFacts}

\begin{defa}[rank] Given an element $\sigma \in R_n$, the {\em rank} of $\sigma$, denoted \textup{rk}$(\sigma)$, is defined to be \textup{rk}$(\sigma) = |\dom(\sigma)| = |\ran(\sigma)|$. 
\end{defa}
There are two main types of elements of $R_n$: the elements of rank $n$ (i.e., the permutations) and the elements of rank less than $n$. In his analysis of the representation theory of $R_n$, Munn \cite{Munn2} introduced what he called {\em cycle-link notation} for the elements of $R_n$. As an example, consider the element $\sigma \in {R_4}$ 
\begin{displaymath}
\sigma = \left( \begin{array}{cccc}
1&2&3&4\\
3&-&2&4 \end{array} \right).
\end{displaymath}
In cycle-link notation, a cycle $(a_1,a_2,\ldots,a_k)$ means that
\[a_1 \mapsto a_2, a_2 \mapsto a_3, \ldots, a_{k-1} \mapsto a_k, \textup{ and } a_k \mapsto a_1,\]
and a link $[b_1, b_2, \ldots, b_k]$ means that 
\[b_1 \mapsto b_2, b_2 \mapsto b_3, \ldots, b_{k-1} \mapsto b_k, \textup{ and } b_k \textup{ goes nowhere.}\]
The element $\sigma$, expressed in cycle-link notation, would be $[1,3,2](4)$. Note that the cycle-link representation of a permutation always consists of cycles only, and elements of rank less than $n$ always contain links when written in cycle-link notation. 

\begin{thm}
\[|R_n| = \sum_{k=0}^{n} {\binom{n}{k}}^2 k!\]
\end{thm}
\begin{proof}
For any particular rank $k$, there are $\binom{n}{k}$ choices for the domain and $\binom{n}{k}$ choices for the range of an element of $R_n$, and for any particular choice of domain and range, there are $k!$ ways of mapping the domain to the range.\end{proof}
We also have the recursive formula
\begin{thm}
\label{SizeRnRecursiveThm}
For $n \geq 3$, 
\begin{equation}
|R_n| = 2n|R_{n-1}| - (n-1)^2|R_{n-2}|.   \label{SizeRnRecursive}
\end{equation}
\end{thm}
\begin{proof}See Section \ref{FactorizationProof}.\end{proof}

\begin{thm}
$R_n$ is an inverse semigroup.
\end{thm}
\begin{proof}
For an element $\sigma \in R_n$, define $\gamma \in R_n$ by 
\[\textup{dom}(\gamma) = \textup{ran}(\sigma), \textup{ and, for } x \in \textup{dom}(\gamma), \]
\[\gamma(x) = \sigma^{-1}(x).\]
It is easy to see that $\gamma\sigma\gamma = \gamma$ and $\sigma\gamma\sigma = \sigma$, and that $\gamma$ is the only element of $R_n$ satisfying both equations.
\end{proof}

\subsection{Semigroup Representations}
\label{SemigpReps}
Here we provide a brief introduction to the representation theory of semigroups. For a more complete introduction to this topic, see \cite{Rhodes} and \cite{Steinberg2}. 

In this paper, we adopt the convention that all maps act on the left of sets, and we only consider left modules. We also only concern ourselves with finite semigroups. For the rest of the paper, every semigroup discussed is understood to be finite. Let $S$ be a semigroup. If $S$ does not have identity, then let $S^1$ denote $S$ with an identity element added. If $S$ does have identity, then let $S^1=S$.

\begin{defa}[semigroup representation] A {\em representation} of $S$ is a finite dimensional $\mathbb{C}$-vector space $V$ together with a homomorphism $\rho$ from $S$ into the semigroup \textup{End}$(V)$ under composition. \end{defa}

\noindent {\bf Remark}: If $S$ is a group, this definition does not agree with the usual group definition, which requires $\rho(e) = 1 \in \textup{End}(V)$. Such representations are called {\em unital}. However, there is only a trivial difference between a representation and a unital representation of $S^1$: a representation is called {\em null} if $\rho(s)=0$ for all $s\in S$, and every representation of $S^1$ is either unital, null, or a direct sum of a unital and a null representation (\cite{Rhodes}, Fact 1.10). 

\begin{defa}[semigroup matrix representation] A {\em matrix representation} of $S$ is a representation of $S$ with respect to a particular choice of basis for $V$, i.e., it is a homomorphism from $S$ into the semigroup of $\dim(V) \times \dim(V)$ matrices with entries in $\mathbb{C}$ under multiplication. 
\end{defa}

\begin{defa}[representation module] A {\em representation module} for $S$ is a $\mathbb CS$-module that is also a finite dimensional $\mathbb C$-vector space.
\end{defa}

\begin{defa}[matrix representation of $\C S$] A {\em matrix representation} of the semigroup algebra $\C S$ is an algebra homomorphism $\rho$ from $\C S$ into the matrix algebra $M_{d_\rho}(\C$). 
\end{defa}

Given a representation $\rho$ of $S$, we may create a representation module for $S$ by viewing the underlying vector space $V$ as a $\mathbb{C}S$-module, where the action of $S$ on $V$ is given by $\rho$, and the action of $\mathbb{C}S$ on $V$ is given by the linear extension of $\rho$. In this case, we call $V$ (as a $\mathbb{C}S$-module) the {\em representation module for $\rho$}. Conversely, given a $\mathbb{C}S$-module $V$ which is also a finite dimensional $\mathbb{C}$-vector space, we may define a representation $\rho : S \rightarrow$ End($V$) by $\rho(s) \cdot v$ = $s \cdot v$ for all $s \in S$ and $v \in V$. In this case, we call $\rho$ the {\em representation associated to $V$}. Furthermore, a matrix representation of $\C S$ defines a matrix representation of $S$ by restricting to a basis, and any matrix representation of $S$ extends linearly to a matrix representation of $\C S$. Matrix representations of $S$ and of $\C S$ are therefore in one-to-one correspondence, and the notions of representations, matrix representations, and representation modules for $S$ are equivalent. We shall move between these notions freely as needed.

Also, we will often abuse the language and simply refer to a representation or matrix representation $\rho$ of $S$, where the underlying vector space $V$ is understood.

\begin{defa}[irreducible representation] A non-null representation $\rho$ of $S$ is said to be {\em irreducible} if the representation module for $\rho$ is non-null and simple. In other words, with respect to no basis does $\rho$ have the form 
\begin{displaymath}
\rho = 
\left( \begin{array}{cc}
\phi_1 & 0 \cr 
g & \phi_2
\end{array} \right)
\end{displaymath}
for some representations $\phi_1$, $\phi_2$ and some matrix-valued function $g$.
\end{defa}

A matrix representation of $\C S$ is irreducible if and only if the corresponding matrix representation of $S$ is irreducible. 

\begin{defa}[equivalence of representations] Given two matrix representations $\rho_1$ and $\rho_2$ of $S$ \LP or of $\C S$\RP, $\rho_1$ is said to be {\em equivalent} to $\rho_2$ if there exists an invertible matrix $A$ such that $A \rho_1(s) A^{-1} = \rho_2(s)$ for all $s\in S$. Equivalently, two representations of $S$ are equivalent if their representation modules are isomorphic \LP as $\mathbb{C}S$-modules\RP.\end{defa} 

\begin{thm}[decomposing representations]For any $S$ such that ${\mathbb{C}S}$ is semi-simple, any representation of $S$ \LP or of $\C S$\RP\textup{ }is equivalent to a direct sum of irreducible and null representations of $S$ \LP resp. $\C S$\RP. Furthermore, there are only finitely many inequivalent, irreducible representations of $S$ \LP resp. $\C S$\RP.\end{thm}

\begin{proof}
This is \cite{Rhodes}, Theorem 1.18 and Proposition 1.19.
\end{proof}

\begin{thm}[Wedderburn-Artin]
\label{Wedderburn}
Let $\mathbb{C}S$ be semisimple, and let ${\mathcal Y}$ be a complete set of inequivalent, irreducible matrix representations for $S$. For each $\rho \in {\mathcal Y}$, \mbox{$\rho : S \rightarrow M_{d_\rho}(\mathbb C)$}; extend $\rho$ linearly to $\tilde \rho: \mathbb CS \rightarrow M_{d_\rho}(\mathbb C)$. Then the family $\{\tilde \rho\}$ defines an isomorphism of algebras
\begin{equation}\label{theisom1}\mathbb{C}S \rightarrow \bigoplus_{\rho \in {\mathcal Y}}M_{d_\rho}(\mathbb C).\end{equation}
Explicitly,
\[
\sum_{s \in S}f(s)s \mapsto
\bigoplus_{\rho \in {\mathcal Y}} \tilde{\rho}\left(\sum_{s\in S}f(s)s\right) =
\bigoplus_{\rho \in {\mathcal Y}}\left( \sum_{s \in S}f(s)\rho(s)\right).
\]
\end{thm}
\begin{proof}See, for example, \cite{Curtis}, \cite{Rhodes}, \cite{Serre}, or \cite{DennisFarb}, p. 50, exercise 18.\end{proof}

\begin{cor}[sum of the squares of the dimensions] Let $S$ be a semigroup such that $\mathbb CS$ is semisimple, and let ${\mathcal Y}$ be a complete set of inequivalent, irreducible matrix representations for $S$. Then
\begin{equation}
|S| = \sum_{\rho\in\mathcal{Y}} d_\rho^2  
\label{Sum of Squares of Dimensions}
\end{equation}
\end{cor}
\begin{proof}The formula (\ref{Sum of Squares of Dimensions}) is just the $\mathbb C$-dimensionality of the algebras appearing in (\ref{theisom1}).\end{proof}

\begin{thm}[Munn \cite{Munn3}, Theorem 3.1]$\mathbb CR_n$ is semisimple.\end{thm}
More generally, we have the following extension of Maschke's theorem. 
\begin{thm}[Munn \cite{Munn2}, Theorem 4.4]If $S$ is an inverse semigroup, then $\mathbb CS$ is semisimple.\end{thm}


\subsection{Adapted Representations}
\label{AdaptedReps}
Let $X_0 < X_1 < \ldots < X_n$ be a chain of semigroups whose semigroup algebras $\C X_i$ are semisimple. Even if a matrix representation of $X_i$ is irreducible, its restriction to $X_{i-1}$ will typically not be. However, it will be equivalent to (though not necessarily \emph{equal} to) a direct sum of irreducible and null matrix representations of $X_{i-1}$. 

\begin{defa}[adapted representations] Let ${\Y}_i$ be a set of inequivalent, irreducible matrix representations for $X_i$. The collection $\{{\mathcal Y}_i\}_{i=0}^n$ is said to be {\em chain-adapted to the chain} $X_0<X_1<\ldots <X_n$ if, for every $i \geq 1$ and every $\rho \in {\Y}_i$, $\rho|_{X_{i-1}}$ is \emph{equal} to a direct sum of null representations and matrix representations in ${\mathcal Y}_{i-1}$. Note that this definition forces the irreducible representations appearing in such a restriction to be equal whenever they are equivalent. \end{defa}

Induction shows that a representation from a chain-adapted set of representations may be restricted further down the chain with the same equality results. In particular, if $\{{\mathcal Y}_i\}_{i=0}^n$ is chain-adapted to $X_0 < X_1 < \ldots < X_n$, then $\{{\mathcal Y}_j\}_{j\in J}$ is chain adapted to $<_{j \in J}X_j$ for any nonempty subset $J \subseteq \{1, \ldots, n\}.$

We shall often simply ask that a complete set of irreducible representations ${\mathcal Y}_n$ for $X_n$ be adapted to the chain $X_0 < X_1 < \ldots < X_n$. In this case, the choices of the ${\mathcal Y}_i$ are understood (and, in fact, are often completely determined) by the choice of ${\mathcal Y}_n$.

\bigskip \noindent {\bf Remark}: Adapted representations are very important in the construction of FFTs (see \cite{DanDiameters}, for example), but the requirement that a set of representations be adapted is in no way limiting. Given any chain $X_0 < X_1 < \ldots < X_n$ of semigroups whose semigroup algebras $\mathbb C X_i$ are semisimple, a straightforward induction argument shows that a complete set of inequivalent, irreducible, chain-adapted matrix representations always exists. Explicitly describing such sets, however, is frequently a challenging endeavor.

\bigskip \noindent {\bf Remark}: Chain-adapted representations are also sometimes referred to as {\em seminormal} or {\em Gelfand-Tsetlin} representations.

\subsection{Schur's Lemma}

\label{Schurs}

Schur's Lemma appears in many forms, the most widely recognized probably being the one below.

\begin{lem}[Schur's Lemma]
\label{Schurs1}
Let $R$ be a ring and let $M, N$ be simple $R$-modules. Then every $R$-module homomorphism $\phi : M \rightarrow N$ is either 0 or an isomorphism. 
See, e.g. \cite{DennisFarb}, p. 30-32.
\end{lem}

For computational purposes, the following (also called Schur's Lemma) is a very useful consequence of Lemma \ref{Schurs1}.


\begin{lem}[Schur's Lemma]
Suppose that $A < B < C$ are semigroups, that ${\mathcal Y}_A$, ${\mathcal Y}_B$, and ${\mathcal Y}_C$ are complete sets of inequivalent, irreducible matrix representations for $A, B, C$ respectively, adapted to the chain, and that $\C A, \C B$, and $\C C$ are semisimple. Let $\rho \in {\mathcal Y}_C$. 
Say $\rho|_B = \rho^1 \oplus \cdots \oplus \rho^k$ \LP with each $\rho^j$ either in ${\Y}_B$ or, without loss of generality, null of dimension 1\RP\textup{ }and $\rho^j|_A = \rho^j_1 \oplus \cdots \oplus \rho^j_{g(j)}$ \LP with each $\rho^j_i$ either in ${\Y}_A$ or, again without loss of generality, null of dimension 1\RP. Let $\sigma \in B$ such that $\sigma$ commutes with A. Then $\rho(\sigma)$ is a block matrix
\[
\rho(\sigma) = \left( 
\begin{array}{ccccc}
\rho^1(\sigma) & 0 & 0 & \ldots & 0 \cr
0 & \rho^2(\sigma) & 0 & \ldots & 0 \cr
0 & 0 & \ddots & \ddots & \vdots \cr 
\vdots & \vdots & \ddots & \ddots & 0 \cr
0 & 0 & \ldots & 0 & \rho^k(\sigma) 
\end{array} 
\right)
\]
where, for $1\leq j\leq k$, the block $\rho^j(\sigma)$ is itself a block matrix, with blocks indexed by all ordered pairs from $\{\rho_1^j, \ldots , \rho_{g(j)}^j\}$:
\begin{displaymath}
\rho^j(\sigma) = \bordermatrix{&\rho^j_1&\rho^j_2&\ldots&\rho^j_{g(j)} \cr
          \rho^j_1 & \lambda_{1,1}^j I  & \lambda_{1,2}^j I  & \ldots  & \lambda_{1,g(j)}^j I  \cr
          \rho^j_2 & \lambda_{2,1}^j I  & \lambda_{2,2}^j  & \ldots  & \lambda_{2,g(j)}^j I  \cr
          \vdots & \vdots  & \vdots  & \ddots  & \vdots  \cr
          \rho^j_{g(j)} & \lambda_{g(j),1}^j I  & \lambda_{g(j),2}^j I  & \ldots  & \lambda_{g(j),g(j)}^j I \cr}
\end{displaymath}
where the $\lambda_{a,b}^j$ are scalars, $I$ is the appropriately-sized identity matrix for the block in which it appears, and the block in position $\rho_{a}^j, \rho_{b}^j$ is non-zero only if $\rho_{a}^j$ and $\rho_{b}^j$ are equivalent representations of $A$.
\end{lem}

\begin{proof}
To prove this, we need only show that $\rho^j(\sigma)$ has the block form indicated above. Let $V^j$ be the representation module for $\rho^j$ on $B$. Viewing $V^j$ as a $\C A$-module, we have $V^j = V_1^j \oplus \cdots \oplus V_{g(j)}^j$, and the $V_i^j$ are simple $\C A$-modules where the module action on $V_i^j$ is given by $\rho_i^j$. Note that a basis $\mathcal B$ of $V^j$ has already been chosen by virtue of the fact that $\rho^j$ is in matrix form, and likewise, bases ${\mathcal B}_i$ for the $V_i^j$ are given since $\rho_i^j$ is in matrix form. Since $\rho^j$ is adapted, we have that ${\mathcal B}_i \subseteq {\mathcal B}$ for all $i$, ${\mathcal B}_i \cap {\mathcal B}_k = \emptyset$ for $i \neq k$, the ${\mathcal B}_i$ are ordered as subsets of ${\mathcal B}$, and $\cup {\mathcal B}_i = {\mathcal B}$. 

Since $\sigma$ commutes with A, we have that, for all $v \in V^j$ and $x \in A$, 
\[
\rho^j(\sigma) \cdot x \cdot v = \rho^j(\sigma x) \cdot v = \rho^j(x\sigma) \cdot v = x \cdot \rho^j(\sigma) \cdot v,
\]
i.e., $\rho^j(\sigma)$ is a $\C A$-linear map from $V^j$ to itself. Hence 
\[
\rho^j(\sigma) \in \textup{Hom}_{\C A}(V^j, V^j) = \textup{Hom}_{\C A}(V_1^j \oplus \cdots \oplus V_{g(j)}^j, V_1^j \oplus \cdots \oplus V_{g(j)}^j)
.\]
According to our basis ${\mathcal B}$ for $V^j$, we have
\begin{small}
\[
\rho^j(\sigma) \in 
\left( 
\begin{array}{cccc}\textup{Hom}_{\C A}(V^j_1,V^j_1)  &   \textup{Hom}_{\C A}(V^j_2,V^j_1)  &  \ldots  &  \textup{Hom}_{\C A}(V^j_{g(j)},V^j_1)  \cr
\textup{Hom}_{\C A}(V^j_1,V^j_2)  &   \textup{Hom}_{\C A}(V^j_2,V^j_2)  &  \ldots  &  \textup{Hom}_{\C A}(V^j_{g(j)},V^j_2)  \cr
\vdots & \vdots & \ddots & \vdots \cr
\textup{Hom}_{\C A}(V^j_1,V^j_{g(j)})  &   \textup{Hom}_{\C A}(V^j_2,V^j_{g(j)})  &  \ldots  &  \textup{Hom}_{\C A}(V^j_{g(j)},V^j_{g(j)})
\end{array}\right).
\]
\end{small}

Let $X \in \textup{Hom}_{\C A}(V^j_a,V^j_b)$ be given in matrix form with respect to the bases ${\mathcal B}_a$, ${\mathcal B}_b$. By Lemma \ref{Schurs1}, $X$ is either $0$ or an isomorphism (i.e. $X$ is either $0$ or $\rho^j_a$ and $\rho^j_b$ are equivalent representations of $A$). Suppose then that $X$ is an isomorphism. The goal, then, is to show that $X$ is a diagonal matrix. Since $X$ is $\C A$-linear, for every $x \in \C A, v\in V^j_a$ we have
\[
Xx\cdot v = x \cdot Xv.
\]
With respect to ${\mathcal B}_a$, ${\mathcal B}_b$ this is 
\[
X \rho^j_a(x) v = \rho^j_b(x) X v,
\]
and hence for every $x \in \C A$
\[
X \rho^j_a(x) = \rho^j_b(x) X.
\]
Since $\rho^j_a$ and $\rho^j_b$ are equivalent, and are either null or part of an adapted set of matrix representations, we have $\rho^j_a = \rho^j_b$, and thus
\begin{equation}
\label{CenterOfMnC}
X \rho^j_a(x) = \rho^j_a(x) X
\end{equation}
for all $x \in \C A$.

Now, either $\rho^j_a$ is null of dimension $1$ (in which case so is $\rho^j_b$), which means that $X$ is $1$-dimensional and we're done, or $\rho^j_a$ is an irreducible representation of $A$. So, suppose $\rho^j_a$ is irreducible. Then by Burnside's theorem (Theorem 1.14 of \cite{Rhodes}), we have
\[
\rho^j_a(\C A) = M_{|{\mathcal B}_a|}(\C),
\]
and therefore (\ref{CenterOfMnC}) says that $X$ is in the center of $M_{|{\mathcal B}_a|}(\C)$, i.e. $X$ is diagonal.
\end{proof}

\noindent {\bf Remark/Notation}: Schur's Lemma says that adapted representations cause the matrix $\rho(\sigma)$ to be sparse and structured. Under the same hypotheses as above, let ${\mathcal{M}}(B,A)$ be the maximum multiplicity of an irreducible or dimension $1$-null representation of $A$ occurring in the restriction, from $B$ to $A$, of an irreducible representation of $B$. Since $\rho(\sigma)$ is a block matrix of the form
\[
\rho(\sigma) = \left( 
\begin{array}{ccccc}
\rho^1(\sigma) & 0 & 0 & \ldots & 0 \cr
0 & \rho^2(\sigma) & 0 & \ldots & 0 \cr
0 & 0 & \ddots & \ddots & \vdots \cr 
\vdots & \vdots & \ddots & \ddots & 0 \cr
0 & 0 & \ldots & 0 & \rho^k(\sigma) 
\end{array} 
\right),
\]
and the block $\rho^j(\sigma)$ is of the form
\begin{displaymath}
\rho^j(\sigma) = \bordermatrix{&\rho^j_1&\rho^j_2&\ldots&\rho^j_{g(j)} \cr
          \rho^j_1 & \lambda_{1,1}^j I  & \lambda_{1,2}^j I  & \ldots  & \lambda_{1,g(j)}^j I  \cr
          \rho^j_2 & \lambda_{2,1}^j I  & \lambda_{2,2}^j  & \ldots  & \lambda_{2,g(j)}^j I  \cr
          \vdots & \vdots  & \vdots  & \ddots  & \vdots  \cr
          \rho^j_{g(j)} & \lambda_{g(j),1}^j I  & \lambda_{g(j),2}^j I  & \ldots  & \lambda_{g(j),g(j)}^j I \cr}
\end{displaymath}
where $\lambda^j_{a,b} \neq 0$ implies $\rho_a^j$ and $\rho_b^j$ are equivalent representations of $A$, we see that the matrix $\rho(\sigma)$ contains at most ${\mathcal{M}}(B,A)$ non-zero entries per row and column. The computational implication of Schur's Lemma, then, is that for an arbitrary $d_\rho \times d_\rho$ matrix $H$, it requires no more than ${\mathcal M}(B,A)d_\rho^2$ complex multiplications and additions to perform each of the matrix multiplications $\rho(\sigma)H$ and $H \rho(\sigma)$, as opposed to the upper bound of $d_\rho^3$ multiplications and additions that would be necessary if $\rho(\sigma)$ were arbitrary.


\section{Basis considerations for $\C S$}
\label{BasisConsid}

\subsection{The poset structure of an inverse semigroup}
\label{Poset}

\begin{defa}[poset structure of $S$] Let $S$ be a finite inverse semigroup. For $s, t \in S$, define
\begin{align*}s \leq t &\iff s = et \textup{ for some idempotent } e\in S \cr 
& \iff s = tf \textup{ for some idempotent } f\in S.
\end{align*}
\end{defa}

For $R_n$, the idempotents are the restrictions of the identity map. If $S$ is a group, then its poset structure is trivial. 

If $P$ is a finite poset, then the {\em zeta function} $\zeta$ of $P$ is given by
\[
\zeta:P\times P \rightarrow \{0,1\}
\]
\[ \zeta(x,y) =
\begin{cases}
1 & \textup{ if }x \leq y \\
0 & \textup{ otherwise.}\end{cases}
\]

Given a poset $P$, one may define an {\em incidence algebra} for $P$ over any ring with identity. The element $\zeta$ is invertible in the incidence algebra, and its inverse is called the {\em M\"obius function} $\mu$. There is a general theory of M\"obius inversion for incidence algebras. We will not go into the details here. Rather, we will record only the results that we need. Details may be found in \cite{Steinberg2}. 

The M\"obius function for $R_n$ (over $\C$) is well known \cite{Stanley}, \cite{Steinberg2}. It is

\[\mu(x,y) = (-1)^{\rk(x)-\rk(y)}.\]

\subsection{Natural bases for $\C S$}
\label{NaturalBases}

Let $S$ be an inverse semigroup. There are two natural bases for $\C S$ (if $S$ is not a group). The first basis is of course $\{ s \}_{s\in S}$, and multiplication in $\C S$ with respect to this basis is just the linear extension of the multiplication in $S$. To motivate the second basis, recall that every (finite) inverse semigroup is isomorphic to a sub-semigroup of a rook monoid and can therefore be viewed as a collection of partial functions. There is another model for composing partial functions: only allow the composition if the range of the first function ``lines up" with the domain of the second. For example, if
\[
\sigma=\left(
\begin{array}{cccc}
1&2&3&4\\
2&-&1&-\\
\end{array}
\right),
\qquad
\pi=\left(
\begin{array}{cccc}
1&2&3&4\\
4&3&-&-\\
\end{array}
\right),
\]
then the idea is that the composition $\pi \circ \sigma$ is
\[
\pi \circ \sigma =
\left(
\begin{array}{cccc}
1&2&3&4\\
3&-&4&-\\
\end{array}
\right),
\]
and the composition $\sigma \circ \pi$ is disallowed. The {\em groupoid basis} for $\C S$ encodes this.

\begin{defa}[groupoid basis] Let $S$ be an inverse semigroup. Define, for each $s \in S$, the element $\ld s \rd \in \mathbb C S$ by
\[\ld s \rd = \sum_{t\in S: t \leq s} \mu(t,s) t.\]
\end{defa}

\begin{thm} The collection $\{\ld s \rd \}_{s\in S}$ is a basis for $\mathbb C S$. Multiplication in $\mathbb C S$ relative to this basis is given by the linear extension of
\begin{equation}
\ld s \rd \ld t \rd = 
\begin{cases}
\ld st \rd & \textup{ if }\dom (s) = \ran (t) \\
0 & \textup { otherwise.}
\end{cases}
\label{GroupoidBasisMult}
\end{equation}
Furthermore, the change of basis to the $\{s\}_{s\in S}$ basis of $\C S$ is given by M\"obius inversion:
\begin{equation}
\label{sissumofbrackett}
s = \sum_{t \in S: t \leq s} \ld t \rd.
\end{equation}
\end{thm}

\begin{proof}
This is \cite{Steinberg2}, Lemma 4.1 and Theorem 4.2, using our convention that maps act on the left of sets.
\end{proof}

The viewpoint then is that we have two natural bases for $\mathbb C S$, the basis $\{s\}_{s \in S}$, and the basis $\{\ld s\rd \}_{s \in S}$. Note that, if $S$ is a group, then $s = \s \in \C S$ for all $s\in S$.


\subsection{Fourier bases for $\C S$}
\label{FourierBases}

Let $S$ be a semigroup such that $\mathbb C S$ is semisimple. This means that $\mathbb C S$ decomposes into a direct sum of simple submodules (i.e. left ideals):
\[
\mathbb C S = \bigoplus \mathbb C L_i.
\]

Let ${\mathcal Y}$ be a complete set of inequivalent, irreducible matrix representations of $\C S$. Then, according to Theorem \ref{Wedderburn},
\begin{equation}
\bigoplus_{\rho \in {\mathcal Y}} \rho : \C S \rightarrow
\bigoplus_{\rho \in {\mathcal Y}} M_{d_\rho} (\C)
\label{FourierIsom}
\end{equation}
is an isomorphism of algebras.

There is a natural basis for the algebra on the right: the set of matrices in the algebra with the property that exactly one entry is 1 (the rest being 0). The inverse image of this set is a basis for $\mathbb C S$ called the {\em dual matrix coefficient basis for ${\mathcal Y}$}, or the {\em Fourier basis for $\CS$ according to ${\mathcal Y}$}. When we refer to a Fourier basis for $\CS$, we mean any basis of $\CS$ that can arise in this manner by choosing an appropriate $\Y$. Note that there is a unique Fourier basis for $\CS$ (up to ordering) if and only if every irreducible representation of $S$ has dimension $1$ (as is the case for $S = \ZnZ$, for which the isomorphism (\ref{FourierIsom}) is the usual discrete Fourier transform and the associated Fourier basis is the usual basis of exponential functions).

Consider the natural basis for the algebra on the right. The preimage of a single column of these elements from the $\rho^{th}$ block is a basis $B$ for a submodule of $\mathbb C S$, and each element of $\mathbb C S$ acts on $B$ exactly as described by $\rho$. We therefore have that $B$ is a basis for an irreducible submodule of $\mathbb C S$ (isomorphic to the representation module for $\rho$). Since the map above was an isomorphism, the preimages of distinct columns have intersection $\{0\}$, and we therefore have the well-known fact: 

\begin{thm}Each irreducible submodule of $\CS$ occurs in the decomposition of $\CS$ into irreducibles exactly as many times as its dimension.
\end{thm}



\section{The Fourier transform on $\C S$}
\label{FourierTransformCS}

Let $S$ be an inverse semigroup, and let $\Y$ be any complete set of inequivalent, irreducible matrix representations of $\C S$.

\begin{defa}The isomorphism \LP \textup{\ref{FourierIsom}}\RP\textup{ }is called a {\em Fourier transform} on $\C S$. 
\end{defa}

\begin{defa}Let $f\in \C S$. The {\em Fourier transform} of $f$ is the image of $f$ in the matrix algebra 
\[\bigoplus_{\rho \in {\mathcal Y}}M_{d_\rho} (\C)\]
via the isomorphism \LP\textup{\ref{FourierIsom}}\RP. Equivalently, the Fourier transform of $f$ is the re\hyph expression of $f$ in $\C S$ in terms of a Fourier basis for $\mathbb C S$.
\end{defa}

Let $f \in \C S$ be given with respect to one of the natural bases, i.e., either
$f = \sum_{s\in S}f(s)s$ or $f=\sum_{s\in S}f(s)\ld s \rd$.
We shall sometimes say ``calculating the Fourier transform on $\C S$" to mean calculating the Fourier transform of an arbitrary element $f \in \C S$ given with respect to one of the natural bases, where the choice of ${\mathcal Y}$ is understood.

\begin{defa}Let $\rho$ be a representation of $\C S$. Define
\[
\hat f(\rho) =
\begin{cases}
\sum_{s\in S}f(s)\rho(s) & \textup{ if } f = \sum_{s\in S}f(s)s, \\
\sum_{s\in S}f(s)\rho(\s) & \textup{ if } f = \sum_{s\in S}f(s)\s. \\
\end{cases}
\]
For $\rho \in \Y$, $\hat f (\rho)$ is therefore just the $\rho^{th}$ block in the image of $f$ in the isomorphism \LP\textup{\ref{FourierIsom}}\RP.
\end{defa}

\noindent {\bf Remark}: Since the map in (\ref{FourierIsom}) is an isomorphism, it respects multiplication, and hence the Fourier transform turns convolution of elements of $\C S$ (\ref{CSConvolutionSemigpBasis}) into multiplication of block-diagonal matrices. It turns convolution into pointwise multiplication if and only if all irreducible representations of $S$ have degree one, which, for example, is the case for $S$ = $\ZnZ$.


\section{Decomposition of $\C R_n$ into a matrix algebra over group algebras}
\label{SteinDecomp}

Let $S$ be an inverse semigroup. The theorem in this section is a special case of Theorem 4.6 in Steinberg \cite{Steinberg2}, which provides an explicit isomorphism between $\CS$ and a direct sum of matrix algebras over group algebras. The groups appearing in the group algebras are the maximal subgroups of $S$. The purpose of this section is to explain this isomorphism in the case when $S = R_n$.

This turns out to be a very important isomorphism for the FFT theory for $R_n$. It allows us to easily describe a complete set of inequivalent, irreducible representations for $R_n$ (Section \ref{SteinReps}), it allows us to describe explicit Fourier bases for $\CRn$ (Section \ref{ExplicitFourierBasis}), it allows us to easily state the Fourier inversion theorem for $\C R_n$ (Section \ref{FourierInversion}), and it forms the basis for one of our FFT algorithms (Section \ref{RnFFT1}).


Let $D_k \subseteq \CRn$ be the $\mathbb C$-span of $\{ \s : s \in R_n, \rk(s) = k\}.$
By (\ref{GroupoidBasisMult}), we have that 
$\CRn = \bigoplus_{i=0}^{n} D_k,$ so that the product of two elements in $D_k$ is an element of $D_k$, and the product of an element in $D_i$ with an element in $D_j$ is $0$ if $i \neq j$.

For $k \in \{0,\ldots,n\}$, 
we identify $S_k \subseteq R_n$ with the subgroup of elements of rank $k$ with domain and range equal to $\{1,\ldots,k\}$. We therefore take $S_0$ to be the set consisting of the zero map, and $S_0 \cong S_1$. The following decomposition theorem was implicit in the work of Munn and made explicit in \cite{Steinberg2}, Theorem 4.6.

\begin{thm}
\label{SteinDecompThm}
$D_k \cong M_{\binom{n}{k}} (\C S_k)$, and thus 
$\CRn \cong \bigoplus_{k=0}^{n} M_{\binom{n}{k}} (\C S_k).$
\end{thm}

\begin{proof}
Use the $k$-subsets of $\{1,\ldots,n\}$ to index the rows and columns of \[M_{\binom{n}{k}} (\C S_k)\] in such a way that $\onetok$ is the first $k$-set. If $A$ and $B$ are two $k$-subsets of $\oneton$, let $p_{(A \rightarrow B)}$ denote the unique order preserving bijection from $A$ to $B$. Keeping in mind that our maps act on the left of sets, define a map 
\[
\phi:D_k \rightarrow M_{\binom{n}{k}} (\C S_k)
\]
by defining it on a basis element $\s$:
\[
\phi(\s) = p_{(\ran(s) \rightarrow \onetok)} s p_{(\onetok \rightarrow \dom(s))} E_{\ran(s),\dom(s)},
\]
where $E_{\ran(s),\dom(s)}$ is the $\binom{n}{k} \times \binom{n}{k}$ matrix with a $1$ in the $\ran(s),\dom(s)$ position and $0$ elsewhere. Observe that 
$p_{(\ran(s) \rightarrow \onetok)} s p_{(\onetok \rightarrow \dom(s))} \in S_k$
by construction.

It is easy to show that $\phi$ is an isomorphism and that $\phi^{-1}$ is induced by
\[
g E_{A,B} \mapsto \ld p_{(\onetok \rightarrow A)} g p_{(B \rightarrow \onetok)} \rd
.\]

Note that if $g\in S_k \subseteq R_n$, then
\[
\phi(\ld g \rd) = gE_{1,1}
\]

\end{proof}

\noindent {\bf Notational remark:} For the rest of the paper, if $A$ and $B$ are $k$-subsets of $\{1,\ldots,n\}$, then we will denote the unique order preserving bijection from $A$ to $B$ by $p_{(A \rightarrow B)}$.


\section{Consequences of Theorem \ref{SteinDecompThm}}
\label{SteinDecompConseq}

\subsection{Explicit matrix representations of $\C R_n$}
\label{SteinReps}

Let $G$ be a finite group. The representations of the matrix algebra $M_n(\C G)$ were studied as early as 1942 by A.H. Clifford in the context of Brandt groupoids \cite{Clifford}. In the notation of Section \ref{SteinDecomp}, we have
\[
\CRn \cong \bigoplus_{k=0}^{n} M_{\binom{n}{k}} (\C S_k).
\]
Given an irreducible matrix representation $\rho$ of $S_k$ (or of $\C S_k$), we can ``tensor up" to an irreducible matrix representation $\bar \rho$ of $M_{\binom{n}{k}} (\C S_k)$ and then extend $\bar \rho$ to an irreducible matrix representation of $\CRn$ by declaring it to be $0$ on the other summands. Specifically, for $g \in S_k$,
\begin{equation}
\label{tensorup}
\bar \rho (g E_{A,B}) = E_{A,B} \otimes \rho(g).
\end{equation}
Let $IRR(S_k)$ be any complete set of inequivalent, irreducible matrix representations of $S_k$. For $\rho_1, \rho_2 \in IRR(S_k)$, \cite{Clifford} shows that $\bar{\rho_1}$ and $\bar{\rho_2}$ are equivalent if and only if $\rho_1$ and $\rho_2$ are equivalent, and that all irreducible matrix representations (up to equivalence) of $M_{d_\rho}(\C S_k)$ are obtained in this manner. 

Therefore, the distinct irreducible representations of $\C R_n$ are in one-to-one correspondence with $\biguplus_{k=0}^{n} IRR(S_k)$. Many explicit descriptions of $IRR(S_k)$ are known. Two well-known computationally advantageous ones are Young's seminormal and Young's orthogonal forms. A description of the former may be found in \cite{Clausen}, and a description of the latter in \cite{Dan1990}.

To be completely explicit, suppose $f \in \C R_n$, with 
$f = \sum_{s\in R_n}f(s)\s$,
and let $\rho \in IRR(S_k)$. Then $\bar\rho$ is an irreducible matrix representation of $\C R_n$, and
\begin{align*}
\bar\rho (f) &= \sum_{s\in R_n}f(s)\bar\rho(\s) \\&=
\sum_{\substack{s\in R_n: \\ \rk(s)=k}}f(s)\bar\rho(\s) \\&=
\sum_{\substack{A \subseteq \oneton \\ |A|=k}}\, \sum_{\substack{B \subseteq \oneton \\ |B|=k}}\, \sum_{g\in S_k}
f(p_{(\onetok \rightarrow A)}g p_{(B\rightarrow \onetok)})
\left( E_{A,B}\otimes\rho(g) \right).
\end{align*}

We have thus effectively described the irreducible representations of $\C R_n$ by describing their actions on the $\{\s\}$ basis. We remark that, if desired, one could obtain a description of the irreducible representations of $R_n$ by looking at these representations on the $\{s\}$ basis via (\ref{sissumofbrackett}). We also remark that there are other ways to describe a complete set of inequivalent, irreducible representations for $R_n$ that do not involve decomposing it into a sum of matrix algebras over group algebras, or indeed referencing the $\{ \s \}$ basis at all. For example, see Grood's description in \cite{Grood}, Halverson's description in \cite{Halverson}, or Section \ref{RnRepresentations}.


\subsection{Explicit Fourier bases for $\C R_n$}
\label{ExplicitFourierBasis}

Let $IRR(S_k)$ be a complete set of inequivalent, irreducible matrix representations for $S_k$, for $k \in \zeroton$. For each $ \rho \in IRR(S_k)$, let $\bar \rho$ denote its extension (via the method discussed in Section \ref{SteinReps}) to $\C R_n$. We take $\Y = \bigcup \bar \rho$. 

It is now easy to explicitly describe the Fourier basis for $\C R_n$ according to $\Y$ in terms of the natural ${\{ \s \}}$ basis. That is, if $B \subseteq \C R_n$ is the set of inverse images of the natural basis of $\bigoplus_{\bar \rho \in \Y}M_{d_{\bar \rho}}(\C)$ in the isomorphism 
\begin{equation}
\bigoplus_{\bar\rho \in \Y}\bar\rho:
\C R_n \rightarrow \bigoplus_{\bar\rho \in \Y}M_{d_{\bar \rho}}(\C)
\label{FourierIsomRn}
,\end{equation}
then for each $b \in B$,
\[b = \sum_{s\in R_n}b(s)\s. \]
We will now describe the coefficients $b(s)$.

We begin by assuming we have an explicit description of a Fourier basis for $\C S_k$ for each $k \in \zeroton$. That is, if $C$ is the set of inverse images of the natural basis of the algebra on the right in the isomorphism
\begin{equation}
\bigoplus_{\rho \in IRR(S_k)}\rho:
\C S_k \rightarrow \bigoplus_{\rho \in IRR(S_k)}M_{d_{\rho}}(\C),
\label{FourierIsomSk}
\end{equation}
then, for each $c\in C$,
\[
c = \sum_{x\in S_k}c(x)x,
\] 
and we assume that we know the coefficients $c(x)$. They may be found, for example, by using the standard Fourier inversion theorem for groups (Theorem \ref{FInvGrp}).

Now, fix $\rho \in IRR(S_k)$. Fix $c_{i,j} \in \C S_k$,
\[
c_{i,j} = \sum_{x\in S_k}c_{i,j}(x)x,
\]
to be the inverse image in the isomorphism (\ref{FourierIsomSk}) of the element of 
\[
\bigoplus_{\rho \in IRR(S_k)}M_{d_{\rho}}(\C)
\]
that is $1$ in the $i,j$ position in the $\rho$ block and $0$ elsewhere. $\bar \rho$ is a block matrix whose rows and columns are indexed by the $k$-subsets of $\oneton$, and whose entries are themselves $d_\rho \times d_\rho$ matrices. Call the $m^{th}$ $k$-set $M$ and the $n^{th}$ $k$-set $N$.

\begin{thm}
\label{FourierBasisRnThm}
Let $X$ be a $d_\rho \times d_\rho$ matrix with a $1$ in the $i,j$ position and $0$ elsewhere. The inverse image in the isomorphism \LP\textup{\ref{FourierIsomRn}}\RP\textup{ }of the element of $\bigoplus_{\bar \rho \in \Y}M_{d_{\bar \rho}}(\C)$ that is $E_{m,n}\otimes X$ in the $\bar \rho$ block and $0$ elsewhere is
\[
\ld p_{(\onetok \rightarrow M)} \rd \left(\sum_{x\in S_k}c_{i,j}(x)\ld x \rd \right)
\ld p_{(N \rightarrow \onetok)} \rd
.\]
\end{thm}

\begin{proof}
If $\bar \gamma \in \Y$, $\bar \gamma \neq \bar \rho$, then  
\[
\bar \gamma \left( \ld p_{(\onetok \rightarrow M)} \rd \left(\sum_{x\in S_k}c_{i,j}(x)\ld x \rd \right)
\ld p_{(N \rightarrow \onetok)} \rd \right) = 0
\]
because 
\[
\bar \gamma \left(\sum_{x\in S_k}c_{i,j}(x)\ld x\rd \right)= 0.
\]
On the other hand,
\[
\bar \rho \left(\ld p_{(\onetok \rightarrow M)} \rd \left(\sum_{x\in S_k}c_{i,j}(x)\ld x \rd \right)
\ld p_{(N \rightarrow \onetok)} \rd \right)
\]
\[
=\left(E_{m,1}\otimes I_{d_\rho} \right)
\left(E_{1,1}\otimes X \right)
\left(E_{1,n}\otimes I_{d_\rho} \right)
=E_{m,n}\otimes X.
\]
\end{proof}


\subsection{The Fourier inversion formula for $\C R_n$}
\label{FourierInversion}

We begin by recalling the Fourier inversion theorem for groups.
\begin{thm}
\label{FInvGrp}
Let $G$ be a finite group, and $f=\sum_{s\in G}f(s)s \in \C G$. Let $IRR(G)$ be a complete set of inequivalent, irreducible matrix representations for $G$. Then
\[
f(s) = \frac{1}{|G|} \sum_{\rho \in IRR(G)} d_\rho \textup{trace} \left( \hat f(\rho)\rho(s^{-1}) \right).
\]
\end{thm}

\begin{proof}
See \cite{Serre}, Section 6.2.
\end{proof}

We now state a preliminary lemma, which is a direct consequence of Theorem \ref{FInvGrp} and Section \ref{SteinReps}. 

\begin{lem}
\label{FInvLemma}
Let $f = \sum_{x\in R_n}f(x)\ld x \rd \in \C R_n$, and let $\Y$ be a complete set of inequivalent, irreducible matrix representations for $R_n$, induced by 
\[
\biguplus_{k \in \zeroton} IRR(S_k)
\]
in the manner described in Section \textup{\ref{SteinReps}}. Let $\rk(x) = k$. Write 
\[
x = p_{(\onetok \rightarrow \ran(x))} y p_{(\dom(x) \rightarrow \onetok)}
\]
for a unique $y \in S_k$. 

The ${\ran(x), \dom(x)}$ entry of $\hat f(\bar \rho)$ is a $d_\rho \times d_\rho$ matrix. If we denote it by $\hat f(\bar \rho)_{\ran(x),\dom(x)}$, then we have
\[
f(x) = \frac{1}{|S_k|} \sum_{\rho \in IRR(S_k)} d_\rho \textup{trace}
\left(
\hat f(\bar \rho)_{\ran(x),\dom(x)} \rho(y^{-1})
\right).
\]
\end{lem}

\begin{proof} 
For $\rho \in IRR(S_k)$ we have 
\[
\hat f(\bar \rho) = \sum_{s\in R_n}f(s)\bar\rho(\ld s \rd),
\]
with $\bar \rho(s) = 0$ if $\rk(s) \neq k$. By definition of $\bar \rho$, the Fourier transform $\hat f(\bar \rho)$ may be computed block by block (and thus may be inverted block by block). The ${\ran(x), \dom(x)}$ entry of $\hat f(\bar \rho)$ is determined by the $f(s)$ for which $\ran(s)=\ran(x)$ and $\dom(s)=\dom(x)$, and such $f(s)$ do not affect any other entries of $\hat f(\bar \rho)$. Explicitly, the ${\ran(x), \dom(x)}$ entry of $\hat f(\bar \rho)$ is given by
\[
\sum_{s\in S_k} f(p_{(\onetok \rightarrow \ran(x))} s p_{(\dom(x)\rightarrow \onetok)}) \rho(s).
\]
Let us define a function $f_{\ran(x), \dom(x)}$ on $S_k$ by
\[
f_{\ran(x), \dom(x)}(s) = f(p_{(\onetok \rightarrow \ran(x))} s p_{(\dom(x)\rightarrow \onetok)})
.\]
Then 
\begin{align*}
\sum_{s\in S_k} f(p_{(\onetok \rightarrow \ran(x))} s p_{(\dom(x)\rightarrow \onetok)}) \rho(s) &= 
\sum_{s\in S_k} f_{\ran(x), \dom(x)}(s)\rho(s) \\&
= \hat f_{\ran(x), \dom(x)}(\rho)
.
\end{align*}
The Fourier inversion theorem for groups then applies, and yields: 
\begin{align*}
f(x) &= 
f(p_{(\onetok \rightarrow \ran(x))} y p_{(\dom(x)\rightarrow \onetok)}) \\&= 
\frac{1}{|S_k|} \sum_{\rho \in IRR(S_k)} d_\rho \textup{trace} \left( \hat f_{\ran(x),\dom(x)}(\rho)\rho(y^{-1}) \right),
\end{align*}
and since
\[
\hat f_{\ran(x),\dom(x)}(\rho) =
\hat f(\bar \rho)_{\ran(x),\dom(x)}
,\]
we are done.
\end{proof}

Now, let ${\mathcal S}$ be {\em any} set of inequivalent, irreducible matrix representations for $R_n$. We define some notation before stating the Fourier inversion theorem for $\C R_n$. 

Let $\Y$ be a complete set of inequivalent, irreducible matrix representations for $R_n$ induced by $\biguplus_{k \in \zeroton} IRR(S_k)$ in the manner described in Section \ref{SteinReps}. If $\rho \in IRR(S_k)$, we have the corresponding $\bar \rho \in \Y$, which is equivalent to some representation in $\mathcal S$, which we denote by $\bar{\bar \rho}$.

\begin{thm}[Fourier inversion theorem for $\C R_n$]
Let 
\[
f = \sum_{x\in R_n}f(x)\ld x \rd \in \C R_n.
\]
Let $\rk(x)=k$, and let us denote the semigroup inverse of $x$ by $x^{-1}$. Then
\[
f(x) = \frac{1}{|S_k|} \sum_{\rho \in IRR(S_k)}d_\rho \textup{trace}
\left(
\hat f (\bar{\bar \rho}) \bar{\bar \rho}(\ld x^{-1}\rd)
\right).
\]
\end{thm}

\begin{proof}
Since $\bar{\bar\rho}$ is equivalent to $\bar \rho$, write
\[
\bar \rho = A^{-1} \bar{\bar\rho} A
\]
for some invertible matrix $A$. We therefore have
\[
\hat f(\bar \rho) = A^{-1}\hat f(\bar{\bar\rho}) A.
\]
As in Lemma \ref{FInvLemma}, let $y$ be the unique element of $S_k$ such that
\[
x = p_{(\onetok \rightarrow \ran(x))} y p_{(\dom(x) \rightarrow \onetok)}.
\]
Now, we have
\begin{align*}
\textup{trace}\left(
\hat f(\bar \rho)_{\ran(x),\dom(x)} \rho(y^{-1})
\right)
=&
\textup{trace} \left(\left[
\hat f(\bar\rho) \right] \left[ E_{\dom(x),\ran(x)}\otimes \rho(y^{-1}) \right]
\right) \\ =&
\textup{trace} \left(
\hat f(\bar \rho)\bar\rho(\ld x^{-1}\rd)
\right) \\ =&
\textup{trace} \left(\left[
A^{-1}\hat f(\bar{\bar \rho}) A\right] \left[ A^{-1} \bar{\bar\rho}(\ld x^{-1}\rd)A\right]
\right) \\ =&
\textup{trace} \left(
A^{-1}\hat f(\bar{\bar \rho})\bar{\bar\rho}(\ld x^{-1}\rd)A 
\right) \\ =&
\textup{trace} \left(
\hat f(\bar{\bar \rho})\bar{\bar\rho}(\ld x^{-1}\rd)
\right),
\end{align*}
the last equality arising from the similarity-invariance of trace. The theorem now follows from Lemma \ref{FInvLemma}.
\end{proof}

\section{Associating functions on $S$ to elements of $\C S$}
\label{FunctionAssoc}

We are now ready to address the issues that arise in choosing how to associate functions on $S$ to elements of the semigroup algebra $\C S$. We also define the Fourier transform and the convolution of functions on $S$. 

Let $S$ be an inverse semigroup. Let $f: S \rightarrow \C$. Having defined two ``natural" bases for $\C S$, we have two natural choices for how to associate $f$ to an element of $\C S$, either 
\[
f \leftrightarrow \sum_{s\in S} f(s)s 
\quad
\textup{ or }
\quad
f \leftrightarrow \sum_{s\in S} f(s)\s .
\]

In the first case, which we shall call the {\em semigroup association model}, the elements $s$ in $\C S$ are associated to the characteristic functions of the elements $s \in S$. In the second case, which we shall call the {\em groupoid association model}, the elements $\s$ in $\C S$ are associated to the characteristic functions of the elements $s \in S$. 

\begin{defa}[Fourier transform of a function on $S$] Let $f:S \rightarrow \C$. Given an association model, the Fourier transform of $f$ is defined to be the Fourier transform of the associated element of $\C S$.
\end{defa}

Thus, if the groupoid association model is used, then the Fourier transform of $f$ is
\[
\bigoplus_{\rho \in \Y}\sum_{s\in S}f(s)\rho(\s),
\]
and if the semigroup association model is used, then the Fourier transform of $f$ is
\[
\bigoplus_{\rho \in \Y}\sum_{s\in S}f(s)\rho(s),
\]
where $\Y$ is a complete set of inequivalent, irreducible matrix representations of $\C S$. Convolution of functions on $S$ is also defined by the association model.

\begin{defa}[convolution of functions on $S$] Let $f,g:S \rightarrow \C$. Choosing an association model defines how the images of $f$ and $g$ in $\C S$ multiply. Denote the images of $f, g$ in $\C S$ by $\bar f, \bar g$ respectively. Then we define, for $s\in S$:

Under the semigroup association model, $f\ast g (s) =$ 
the $s^{th}$ coefficient of \mbox{$\bar f \ast \bar g$} expressed with respect to the $\{s\}$ basis.
Under the groupoid association model, \mbox{$f\ast g (s) =$} 
the $\s^{th}$ coefficient of $\bar f \ast \bar g$ expressed with respect to the $\{\s\}$ basis.
\end{defa}

For the case $S = R_n$, what follows is an overview of the considerations one should take into account when choosing which function association model to use.

\begin{consideration}[Convolution]
Convolution of functions on $R_n$ is defined by the association model used. If the semigroup association model is used, then convolution of $\C$-valued functions $f, g$ on $R_n$ is defined by (\ref{CSConvolutionSemigpBasis}). Specifically,
\[
(f \ast g) (s) = \sum_{r,t \in R_n : rt = s} f(r)g(t).
\]

If the groupoid association model is used, then convolution of $\C$-valued functions $f, g$ on $R_n$ is defined by their multiplication in $\C R_n$ with respect to the $\{\s\}$ basis. Specifically,
\[
(f \ast g) (s) = \sum_{\substack{r \in R_n : \\ \ran(s) = \ran(r)}}f(r)g(r^{-1}s).
\]
\end{consideration}

\begin{consideration}[Inner products]
There is a natural inner product for $\C$-valued functions $f,g$ on $R_n$:
\[
<f,g> = \sum_{s \in R_n} f(s)\overline{g(s)}.
\]
It would be nice to use this inner product. Choosing an inner product on functions and an association model induces an inner product on $\C R_n$ (or equivalently, choosing an inner product on $\C R_n$ and an association model induces an inner product on functions on $R_n$). Specifically, if we choose the above inner product, then the semigroup association model induces the inner product $<\cdot,\cdot>_1$ on $\C R_n$ defined by 
\[
<s,t>_1 =
\begin{cases} 
1 & \textup{if }s=t \\
0 & \textup{otherwise,}
\end{cases}
\]
and the groupoid association model induces the inner product $<\cdot,\cdot>_2$ on $\C R_n$ defined by
\[
<\s,\ld t \rd>_2 =
\begin{cases} 
1 & \textup{if }\s= \ld t \rd \\
0 & \textup{otherwise.}
\end{cases}
\]
However, choosing the inner product $<\cdot,\cdot>_1$ leads to undesirable consequences with regards to the Fourier basis for $\C R_n$, as described in Consideration \ref{3}. For this reason, if the semigroup association model is used, then it is recommended to take the inner product $<\cdot,\cdot>_2$ on $\C R_n$ and use that to induce an inner product on $\C$-valued functions on $R_n$.
\end{consideration}

\begin{consideration}[Orthogonality of isotypic subspaces of $\C R_n$]
\label{3}
Since $\C R_n$ is semisimple, we have
\[
\C R_n = \bigoplus L_i,
\]
where each $L_i$ is simple. Let us group them according to their isomorphism classes:
\[
\C R_n = \bigoplus_{\bar \rho \in \Y} V_{\bar \rho},
\]
where $\Y$ is a complete set of inequivalent, irreducible matrix representations for $R_n$, and $V_{\bar \rho}$ is the sum of all simple submodules of $\C R_n$ isomorphic to the representation module for $\bar \rho \in \Y$. The $V_{\bar\rho}$ are called the {\em isotypic subspaces} or {\em isotypic submodules} of $\C R_n$. Note that this decomposition does not depend on the choice of $\Y$. Let $v \in V_{\bar \rho}, v' \in V_{\bar \rho '}.$ We would like to have an inner product $(\cdot,\cdot)$ on $\C R_n$ such that, if $V_{\bar \rho} \neq V_{\bar \rho '}$, then 
\[
(v,v') = 0,
\]
{\em i.e.}, an inner product under which the isotypic subspaces of $\C R_n$ are orthogonal. 
The claim is that $<\cdot,\cdot>_2$ accomplishes this, and that $<\cdot,\cdot>_1$ in general does not.

\begin{thm}Let $V_{\bar \rho} \neq V_{\bar \rho '}$. Then, using the notation above, $<v,v'>_2 = 0.$
\label{IsotypicsOrthogonal}
\end{thm}
\begin{proof}
By linearity, it suffices to show this in the case that $v$ and $v'$ are Fourier basis elements of $\C R_n$. We take $\Y$ to be the set induced by $\biguplus_{k \in \zeroton}IRR(S_k)$ in the manner described in Section \ref{SteinReps}, and we assume that $v$ and $v'$ are part of a Fourier basis for $\C R_n$ according to $\Y$.

We know that
\[
\C S_k = \bigoplus_{\rho \in IRR(S_k)}W_\rho,
\]
where $W_\rho$ contains all the irreducible submodules of $\C S_k$ isomorphic to the representation module for $\rho$. If $w\in W_\rho, w' \in W_{\rho '}, W_\rho \neq W_{\rho '},$ then under the inner product $[\cdot,\cdot]$ on $\C S_k$ defined by
\[
[\sum_{s\in S_k}a(s)s,\sum_{s\in S_k}b(s)s] = \sum_{s\in S_k} a(s)\overline{b(s)},
\]
it follows from the discussion in Chapter 2 of \cite{Serre} that we have $[w,w'] = 0$.

Now, suppose that ${\bar\rho}$ was induced by $\rho \in IRR(S_k)$ and that ${\bar\rho'}$ was induced by $\rho' \in IRR(S_j)$. By Theorem \ref{FourierBasisRnThm}, when written in terms of the natural $\{ \s \}$ basis, $v$ contains nonzero coefficients only for the elements $\s$ where $s$ is of rank $k$, and $v'$ contains nonzero coefficients only for the elements $\s$ where $s$ is of rank $j$. Thus, if $k \neq j$, we have $<v,v'>_2 = 0$. Suppose then that $k=j$. By Theorem \ref{FourierBasisRnThm}, we have
\begin{align*}
v =& \ld p_{(\onetok \rightarrow B)} \rd   \sum_{s \in S_k} v(s)\s    \ld p_{(A \rightarrow \onetok)} \rd, \\
v'=& \ld p_{(\onetok \rightarrow B')} \rd  \sum_{s \in S_k} v'(s)\s   \ld p_{(A' \rightarrow \onetok)} \rd, \\
\end{align*}
for $A$, $A'$, $B$, $B'$ some $k$-subsets of $\oneton$, and
\[
\sum_{s \in S_k} v(s)s \in W_\rho, \, 
\sum_{s \in S_k} v'(s)s \in W_{\rho '}
\] 
some Fourier basis elements for $\C S_k$. 

If $A \neq A'$ or $B \neq B'$, it is apparent that $<v,v'>_2 = 0$, so suppose further that $A=A'$ and $B=B'$.

Now, since $V_{\bar \rho} \neq V_{\bar \rho '}$, we have $\rho \neq \rho'$, and we therefore note that
\[
[\sum_{s \in S_k} v(s)s, \sum_{s \in S_k} v'(s)s] = 0.
\]
For convenience, we temporarily denote $p_{(\onetok \rightarrow B)}$ by $p_B$ and $p_{(A \rightarrow \onetok)}$ by $p_A$. We have
\[
<v,v'>_2 = \sum_{s\in S_k} \sum_{t \in S_k} v(s)\overline{v'(t)} <\ld p_B s p_A \rd , \ld p_B t p_A \rd>_2,
\]
and $\ld p_B s p_A \rd = \ld p_B t p_A \rd$ if and only if $s = t$, so
\[
<v,v'>_2 = \sum_{s \in S_k} v(s) \overline{v'(s)}
=[\sum_{s \in S_k} v(s),\sum_{s \in S_k} v'(s)]
=0.
\]
\end{proof}

\begin{rmk}The isotypic subspaces $V_{\bar\rho}$ need not be orthogonal under the inner product $<\cdot,\cdot>_1$. For example, consider $\C R_1$, which has two nonisomorphic irreducible representations, each of degree 1, and hence a unique Fourier basis. It decomposes as
\[
\C R_1 = \left( \C-\textup{span}(\ld \textup{Id} \rd) \right) \oplus
\left( \C-\textup{span}(\ld 0 \rd) \right)
.\]
We have
\[
<\ld \textup{Id} \rd, \ld 0 \rd>_1 = <\textup{Id} - (0), (0)>_1 = -1
.\]
\end{rmk}

\end{consideration}


\section{FFT algorithm considerations}
\label{FFTConsiderations}

\subsection{Computational complexity}
\label{CompComplexity}

If $G$ is a finite group, then $\{g\}_{g\in G}$ indexes the natural basis of $\C G$. If $S$ is an inverse semigroup, then $\C S$ has two natural bases, the $\{s\}_{s\in S}$ basis and the $\{\s\}_{s\in S}$ basis. We therefore define the following two notions of computational complexity for the Fourier transform on $\C S$.



\begin{defa}[Computational complexity] Let $\Y$ be a complete set of inequivalent, irreducible matrix representations for $\C S$. For an arbitrary element $f \in \C S$ expressed with respect to the $\{s\}$ basis
\[
f = \sum_{s\in S}f(s)s,
\]
the minimal number of operations to compute the Fourier transform of $f$, i.e., to compute
\begin{equation}
\hat f(\rho) = \sum_{s\in S}f(s)\rho(s)
\label{FourierTransformAtRhoSemigp}
\end{equation}
for all $\rho \in \Y$, is denoted by $\semigroupcplx (S)$.

For an arbitrary element $f\in \C S$ expressed with respect to the $\{\s\}$ basis
\[
f = \sum_{s\in S}f(s)\s,
\]
the minimal number of operations to compute the Fourier transform of $f$, i.e., to compute
\begin{equation}
\hat f(\rho) = \sum_{s\in S}f(s)\rho(\s)
\label{FourierTransformAtRhoGroupoid}
\end{equation}
for all $\rho \in \Y$, is denoted by $\groupoidcplx(S)$.

Now, let $\Y$ vary over all complete sets of inequivalent, irreducible representations for $\C S$. We define
\[
{\mathcal C}^s(S) = \minup_{\Y}(\semigroupcplx(S)),
\]
\[
{\mathcal C}^{\s}(S) = \minup_{\Y}(\groupoidcplx(S)), \textup{ and}
\]
\[
{\mathcal C}(S) = \minup({\mathcal C}^s(S),{\mathcal C}^{\s}(S)).
\]
\end{defa}

An operation is defined to be a single complex multiplication followed by a complex addition. For the purposes of our analysis, we assume that all representations in ${\mathcal Y}$ are precomputed and stored in memory.



If $G$ is a group, the $\{s\}$ and the $\{\s\}$ bases of $\C G$ are identical, so we may drop the superscripts on the complexity notation. For example,
\[
\semigroupcplx(G) = \groupoidcplx(G) = {\mathcal T}_{\Y}(G).
\]

Now, let $\Y$ be any complete set of inequivalent, irreducible matrix representations for $\C R_n$. A naive implementation of the Fourier transform on $R_n$, i.e., computing (\ref{FourierTransformAtRhoSemigp}) and (\ref{FourierTransformAtRhoGroupoid}) directly, gives
\[ 
{\mathcal C}(R_n) \leq \sum_{\rho \in \Y}|R_n|d_\rho^2,
\]
which, by (\ref{Sum of Squares of Dimensions}), gives

\begin{thma}${\mathcal{C}}(R_n) \leq |R_n|^2. $
\end{thma}

As $n$ grows, this cost quickly becomes prohibitive. 

For the same reason, we also have ${\mathcal C}(G) \leq |G|^2$ for any group $G$, but
many families of groups enjoy results along the lines of ${\mathcal{C}}(G) = O(|G|\log^c|G|)$. Indeed, such upper bounds remain the goal in group FFT theory. There are also groups $G$ for which there currently exist greatly improved (but not $O(|G|\log^c|G|)$) algorithms, such as matrix groups over finite fields or, more generally, finite groups of Lie type \cite{DanDiameters}. It is conjectured that there are universal constants $c_1, c_2$ such that for any group $G$, ${\mathcal C}(G) \leq c_1 |G| \log^{c_2} |G|$ \cite{DanAMS}. It is currently known that, for example, 
\begin{align*}
&{\mathcal C}(\mathbb Z/n \mathbb Z) = O(n \log n), \\
&{\mathcal C}(G) \leq 8.5 |G|\log |G| \textup{ for any supersolvable group } G, \textup{ and} \\
&{\mathcal C}(S_n) = O(|S_n|\log^2|S_n|) \textup{ (where } S_n \textup{ is the symmetric group on } n\textup{ letters)}.
\end{align*}
See \cite{CooleyTukey}, \cite{Baum}, and \cite{Maslen}, respectively, for these results, or see Section \ref{SnFFT} for a description of an $O(|S_n|\log^3|S_n|)$ FFT on $S_n$.


We now state the main results of this paper.

\begin{thm}
\label{mainthmgroupoid}
$
{\mathcal C}^{\s}(R_n) \leq O(|R_n|\log^2|R_n|).
$
\end{thm}

\begin{thm}
\label{mainthm}
For any $\epsilon > 0, {\mathcal C}^s(R_n) = O(|R_n|^{1+\epsilon})$, and hence
\[
{\mathcal C}(R_n) = O(|R_n|^{1+\epsilon}).
\]
\end{thm}

We prove these results in Sections \ref{RnFFT1groupoid} and \ref{RnFFT1semigp}, respectively.


\subsection{FFTs on groups}
\label{FFTsOnGroups}

In this section, we describe a general method used to construct FFTs on groups. This section is basically an elaboration of the main idea in \cite{DanDiameters}, and is included in the interest of keeping this document self-contained. This method is used in the FFT on $S_n$ given in Section \ref{SnFFT} (which itself is used in the FFT for $R_n$ given in Section \ref{RnFFT1}). We also will generalize this method directly to create an FFT for $R_n$ in Section \ref{RnFFT2}. 



Let $G$ be a group and let $H$ be a subgroup of $G$. $H$ partitions $G$ into cosets $yH$. Let $C \subseteq G$ be a full set of coset representatives for distinct cosets (i.e., $G=\bigcup_{y\in C}yH$, and $y_1 \neq y_2 \in C \Rightarrow y_1H \cap y_2H = \emptyset$). Let $\rho$ be a representation of $G$. Then, since $\rho$ is a homomorphism, we have the factorization
\begin{equation}
\label{groupfactorization}
\hat{f}(\rho) = \sum_{g \in G} f(g)\rho(g) = \sum_{y\in C}\rho(y)\sum_{h\in H}f_y(h)\rho(h),
\end{equation}
where $f_y(h) = f(yh)$. 

Let ${\mathcal Y}_G$ and ${\mathcal Y}_H$ be complete sets of irreducible, inequivalent matrix representations for $G$ and $H$, respectively, adapted to the chain $G > H$. We would like to compute $\hat f(\rho)$ for all $\rho \in {\mathcal Y}_G$. The idea is that if we already knew all the $\hat{f_y}(\gamma)$ for all $\gamma$ in ${\mathcal Y}_H$ and $y \in C$, then we could construct the inner sums in (\ref{groupfactorization}) based on how $\rho$ splits when restricted to $H$. (To do this in practice for a particular group $G$, we would need a theorem that says how $\rho$ splits when restricted to $H$. Such theorems are known as {\em branching theorems}.) Since ${\mathcal Y}_H$, ${\mathcal Y}_G$ are adapted, this construction can be done for free, since (say $\rho |_H = \rho_1 \oplus \cdots \oplus \rho_k$)
\begin{displaymath}
\sum_{h \in H}f_y(h)\rho(h) = \sum_{h \in H}f_y(h)\left[\rho_1(h)\oplus \cdots \oplus \rho_k(h)\right] = 
\hat{f_y}(\rho_1) \oplus \cdots \oplus \hat{f_y}(\rho_k).
\end{displaymath}
Then, to finish computing $\hat{f}(\rho)$, we need to multiply $\frac{|G|}{|H|}$ matrices together and add the results. Denote
$\sum_{h\in H}f_y(h)\rho(h)$
by $A_y(\rho)$. Then by choosing our representations ${\mathcal Y}_G$, ${\mathcal Y}_H$ for $G$ and $H$ to be adapted to the chain $G>H$, we have 

\begin{thm}
\[
{\mathcal T}_{{\mathcal Y}_G}(G) \leq \frac{|G|}{|H|} {\mathcal T}_{{\mathcal Y}_H}(H) + M_G(C),
\]
where $M_G(C)$ is the number of operations needed to compute the sums
\[
\sum_{y \in C}\rho(y)A_y(\rho)
\]
for all $\rho \in {\Y}_G$, given the matrices $A_y(\rho)$.
\end{thm}
$M_G(C)$ can be made small by choosing the coset representatives $y \in C$ in such a way that $\rho(y)$ is a sparse matrix (or, more generally, can be factored into a small number of sparse matrices).

Of course, this method can be used recursively with a subgroup chain and a corresponding collection of chain-adapted representations, as in \cite{DanAMS}, \cite{DanDiameters}. 


\subsection{An FFT on the Symmetric Group}
\label{SnFFT}

In this section, we describe an algorithm to compute the FFT on $S_n$. The FFT algorithm for $R_n$ presented in Section \ref{RnFFT1} explicitly uses FFTs on $S_n$, and ideas from this algorithm are used in the FFT on $R_n$ presented in section \ref{RnFFT2}.

Consider $S_n$ and the subgroup chain $S_n > S_{n-1} > \ldots > S_1 =\{1\}$, where $S_k$ is identified with the subgroup of $S_n$ that fixes the points $k+1, \ldots , n$. Choose a complete, inequivalent, irreducible set of matrix representations for $\C S_n$ adapted to this chain and call it ${\mathcal{Y}}_n$. Two well-known choices are Young's orthogonal and Young's seminormal forms (see, for example, \cite{Dan1990}). Let $t_j$ be the transposition $(j-1,j)$. Let $e$ denote the identity element of $S_n$. We use the following set of coset representatives:
\begin{displaymath}
\{T_i:1\leq i \leq n, T_i = \underbrace{e\ldots e}_i t_{i+1}t_{i+2} \ldots t_n \}.
\end{displaymath}
Thus, for example, $T_1 = t_{2}t_{3} \ldots t_n$, $T_{n-1} = t_n$, and $T_n = e$. We therefore obtain
\begin{equation}
\sum_{s \in S_n}f(s)\rho(s) = \sum_{i=1}^n\rho(T_i)\sum_{s \in S_{n-1}}f_{T_i}(s)\rho(s),     \label{SnDecomp}
\end{equation}
which then implies
\begin{equation}
{\mathcal T}_{{\mathcal Y}_n}(S_n) \leq \frac{S_n}{S_{n-1}} {\mathcal T}_{{\mathcal Y}_{n-1}}(S_{n-1}) + M_{S_n}(\{T_i\}). \label{SnComplexity}
\end{equation}

Once we have calculated all the $\hat{f_{T_i}}$ on $S_{n-1}$, then for any $\rho \in {\mathcal Y}_n$, we can reconstruct (for free, since we are using chain-adapted matrix representations) the inner sums in (\ref{SnDecomp}) based on how $\rho$ splits when restricted to $S_{n-1}$. This splitting is described by the branching theorem for $S_n$ (see, for example, \cite{Dan1990}, p. 304).

We now turn to analyzing the $M_{S_n}(\{T_i\})$ term in (\ref{SnComplexity}). Notice that, for $j >2$, $t_j \in S_j$ and $t_j$ commutes with $S_{j-2}$ (and $t_2 \in S_2$ and $t_2$ commutes with $S_1$). It is easy to derive from the combinatorics of Young tableaux (see, for example, \cite{James}) that the maximum multiplicity occurring in the restriction of any irreducible representation of $S_j$ to $S_{j-2}$ is 2. By Schur's Lemma, then, for any $j$ ($2\leq j\leq n$) and any $\rho \in {\mathcal{Y}}_n$, $\rho(t_j)$ contains at most $2$ nonzero entries per row and column.

Fix $\rho \in {\mathcal{Y}}_n$. Since $\rho(T_i)$ = $\rho(t_{i+1})\rho(t_{i+2}) \ldots \rho(t_{n})$, computing $\rho(T_i)A_{T_i}(\rho)$ for an arbitrary matrix $A_{T_i}(\rho)$ may be accomplished by multiplying $A_{T_i}(\rho)$ on the left by $\rho(t_n)$, multiplying the result on the left by $\rho(t_{n-1})$, multiplying the result of that on the left by $\rho(t_{n-2})$, etc. It therefore takes a maximum of $2(n-i)d_\rho^2$ operations to perform the multiplication $\rho(T_i)A_{T_i}(\rho)$ (keeping in mind that multiplying by the identity matrix may be done for free), and once $\rho(T_i)A_{T_i}(\rho)$ has been computed for all $T_i$, it takes a maximum of $(n-1)d_\rho^2$ operations to add the results to give $\sum_{i=1}^n\rho(T_i)A_{T_i}(\rho)$. Letting $\rho$ vary over ${\mathcal{Y}}_n$, then, implies 
\begin{align*}
M_{S_n}(\{T_i\}) & \leq 
\sum_{\rho \in {\mathcal{Y}}_n} \sum_{i=1}^n 2(n-i) d_\rho^2 + \sum_{\rho \in {\mathcal{Y}}_n}(n-1)d_\rho^2 \\ &=
(n-1)(n)|S_n|+(n-1)|S_n| \\ &=
(n^2-1)|S_n|.
\end{align*}
Putting this together with (\ref{SnComplexity}), we obtain:
\begin{displaymath}
{\mathcal T}_{{\mathcal Y}_n}(S_n) \leq n{\mathcal T}_{{\mathcal Y}_{n-1}}(S_{n-1}) + (n^2-1)|S_n|.
\end{displaymath}
Induction on $n$ then yields
\begin{thm}
${\mathcal T}_{{\mathcal Y}_n}(S_n) \leq \frac{2}{3}n(n+1)^2n!$ 
\end{thm}

This is of order $n!(\log n!)^3 = |S_n|\log^3|S_n|$. The algorithm described here is the heart of Clausen's FFT on $S_n$ \cite{Clausen}.

By a more careful analysis of the matrix multiplications involved, Maslen has obtained an algorithm for the FFT on $S_n$ of complexity $O(|S_n|\log^2|S_n|)$ (\cite{Maslen}, Theorem 1.1):

\begin{thm}
\label{MaslensSnFFT}
Let ${\Y}_n$ denote a complete set of inequivalent, irreducible representations for $S_n$ in Young's orthogonal or seminormal form. Then
\[
{\mathcal T}_{{\Y}_n}(S_n) \leq
\frac{3}{4}n(n-1)|S_n|.
\]
\end{thm}


\section{An algorithm for the FFT on $R_n$}
\label{RnFFT1}

In this section, we present the faster of our two algorithms for the FFT on $R_n$. In the process, we prove Theorems \ref{mainthmgroupoid} and \ref{mainthm}.

\subsection{An FFT for changing from the $\{\s\}$ basis to a Fourier basis}
\label{RnFFT1groupoid}

Let $f \in \C R_n$ be an arbitrary element, given with respect to the $\{\s\}$ basis:
\[
f = \sum_{s\in R_n}f(s)\s.
\]

Let ${\Y}_{k}$ denote Young's seminormal (or orthogonal) matrix representations for $S_k$. For our complete set of inequivalent, irreducible representations for $R_n$, we take the set $\Y$ induced by $\biguplus_{k\in\zeroton} {\Y}_{k}$ in the manner described in Section \ref{SteinReps}. If $\bar\rho \in \Y$ was induced by $\rho \in {\Y}_{k}$, then we have that $\hat f (\bar \rho)$ is an $\binom{n}{k} \times \binom{n}{k}$ matrix whose rows and columns are indexed by the $k$-subsets of $\oneton$, with entries themselves $d_\rho \times d_\rho$ matrices. By Theorem \ref{SteinDecompThm} and (\ref{tensorup}), we know that for $s\in R_n$,
\[
\bar\rho(\s) =
\begin{cases}
0 & \textup{if }\rk(s) \neq k \\
E_{ran(s), dom(s)} \otimes \rho(p_{(\ran(s) \rightarrow \onetok)} s p_{(\onetok \rightarrow \dom(s))}) & \textup{ otherwise.}
\end{cases}
\]

Let $A$ and $B$ be $k$-subsets of $\oneton$. Then the $A,B$ entry of $\hat f (\bar \rho)$ is
\begin{equation}
\label{BlockInGroupoidFT}
\hat f(\bar\rho)_{A,B} = \sum_{s\in S_k} f(p_{(\onetok \rightarrow A} s p_{B \rightarrow \onetok)}) \rho(s).
\end{equation}
If we define a function $f_{A,B}$ on $S_k$ by
\[
f_{A,B}(s) = f(p_{(\onetok \rightarrow A} s p_{B \rightarrow \onetok)}),
\]
then (\ref{BlockInGroupoidFT}) is just $\hat f_{A,B}(\rho)$, a Fourier transform on $S_k$. An obvious algorithm presents itself: for each $k$, run $\binom{n}{k}^2$ FFTs on $S_k$. By Theorem \ref{MaslensSnFFT}, 
\[
{\mathcal T}_{{\Y}_{k}}(S_k) \leq
\frac{3}{4}k(k-1)|S_k|.
\]

\begin{thma}[Theorem \ref{mainthmgroupoid}]
$
{\mathcal C}^{\s}(R_n) \leq O(|R_n|\log^2|R_n|).
$
\end{thma}

\begin{proof}By the above algorithm,
\begin{align*}
\groupoidcplx(R_n) &
\leq \sum_{k=0}^n \binom{n}{k}^2 {\mathcal T}_{{\Y}_{k}}(S_k) \\&
\leq \sum_{k=0}^n \binom{n}{k}^2 \frac{3}{4}k(k-1)k! \\&
\leq \frac{3}{4}n(n-1) \sum_{k=0}^n \binom{n}{k}^2 k! \\&
\leq \frac{3}{4}n(n-1) |R_n|. \\
\end{align*}
Since ${\mathcal C}^{\s}(R_n) \leq \groupoidcplx(R_n)$, $|R_n| \geq n!$, and $n = O(\log (n!))$, we are done.
\end{proof}


\subsection{An FFT for changing from the $\{s\}$ basis to a Fourier basis}
\label{RnFFT1semigp}

Let \mbox{$f \in \C R_n$} be an arbitrary element, given with respect to the $\{ s \}$ basis:
\[
f = \sum_{s\in R_n}f(s)s.
\]

Thanks to the algorithm in Section \ref{RnFFT1groupoid}, a simple algorithm for calculating the Fourier transform of $f$ is evident: use the set $\Y$ of inequivalent, irreducible matrix representations for $R_n$ used in Section \ref{RnFFT1groupoid}, re-express $f$ in terms of the $\{\s\}$ basis, and then run the algorithm in Section \ref{RnFFT1groupoid}. Under this approach, the only algorithmic complexity left to consider is the complexity of changing from the $\{s\}$ basis to the $\{\s\}$ basis. Suppose
\[
f = \sum_{s\in R_n}g(s)\s.
\]
Since
\[
s = \sum_{t\in R_n: t \leq s} \ld t \rd,
\]
basic linear algebra gives 
\[
g(s) = \sum_{x\in R_n: x\geq s} f(x).
\]

Given $s\in R_n$ with $\rk(s)=k$, there are $|R_{n-k}|$ elements $x \in R_n$ satisfying $x\geq s$. A naive implementation of the $\{s\} \rightarrow \{\s\}$ change of basis therefore takes at most
\begin{equation}
\label{1sInZetaMatrix}
\sum_{k=0}^n\binom{n}{k}^2k!|R_{n-k}|
\end{equation}
operations. 

If we consider the matrix that encodes this change of basis (which consists only of $0$'s and $1$'s), then (\ref{1sInZetaMatrix}) is the number of $1$'s in this matrix, counted row by row, since a row corresponding to an element $s\in R_n$ of rank $k$ has $|R_{n-k}|$ $1$'s in it. On the other hand, a column corresponding to an element $s\in R_n$ of rank $k$ has $2^k$ $1$'s in it, so counting column by column, we obtain:
\[
\sum_{k=0}^n\binom{n}{k}^2k!|R_{n-k}| = 
\sum_{k=0}^n\binom{n}{k}^2k!2^k,
\]
and thus the $\{s\} \rightarrow \{\s\}$ change of basis takes at most
\[
\sum_{k=0}^n\binom{n}{k}^2k!2^k \leq
2^n \sum_{k=0}^n\binom{n}{k}^2k! =
2^n |R_n|
\]
operations. 

Putting this together with Theorem \ref{mainthmgroupoid}, we obtain
\begin{thma}[Theorem \ref{mainthm}] Let $\epsilon > 0$. Then
\[
{\mathcal C}^{s}(R_n) \leq 
O(|R_n|^{1+\epsilon}).
\]
\end{thma}

\begin{proof}
The algorithm in this section consists of two steps: change from the $\{s\}$ basis to the $\{\s\}$ basis, and then run the algorithm in Section \ref{RnFFT1groupoid}. Thus
\begin{align*}
{\mathcal C}^{s}(R_n) & \leq \semigroupcplx(R_n) \\&
\leq 2^n |R_n| + \groupoidcplx(R_n) \\&
\leq 2^n |R_n| + \frac{3}{4}n(n-1) |R_n|.
\end{align*}
For large enough $n$, we have $2^n \geq \frac{3}{4}n(n-1)$ and $|R_n|^\epsilon \geq n!^\epsilon \geq 2^n$, so
\[
\semigroupcplx(R_n) \leq O(|R_n|^{1+\epsilon}).
\]
\end{proof}

{\noindent \bf Remark:} We have developed improved algorithms for the $\{s\} \rightarrow \{\s\}$ change of basis, but none yet which are $O(|R_n|\log^c|R_n|)$. Also, in Section \ref{RnFFT2} we exhibit a different algorithm for changing from the $\{s\}$ basis to a Fourier basis, which avoids the $\{\s\}$ basis entirely. However, it too is a $O(|R_n|^{1+\epsilon})$ algorithm.

\section{Another algorithm for the FFT on $R_n$}
\label{RnFFT2}

In this section, we present an algorithm for the FFT on $R_n$ that is quite different from the algorithm presented in Section \ref{RnFFT1}. The algorithm in Section \ref{RnFFT1} involved examining the poset structure of $R_n$, decomposing $\C R_n$ into matrix algebras over group algebras, using matrix representations for $R_n$ induced by computationally advantageous matrix representations of $S_k$, and running FFTs on $S_k$. 

In contrast, the algorithm in this section directly generalizes the ideas in Sections \ref{FFTsOnGroups} and \ref{SnFFT}. It involves finding a semigroup factorization of $R_n$ into ``cosets" and using a chain-adapted set of matrix representations to run recursive calculations on the sub-semigroup chain $R_n > R_{n-1} > \ldots > R_1$. The algorithm in this section computes the Fourier transform of an element in $\C R_n$ expressed with respect to the $\{s\}$ basis directly, without referencing the $\{\s\}$ basis. It is not an 
$O(|R_n|\log^c|R_n|)$
algorithm, so we do not consider the $\{\s\} \rightarrow \{s\}$ change of basis necessary to run it on the $\{\s\}$ basis. Even for an element expressed with respect to the $\{s\}$ basis, the algorithm in Section \ref{RnFFT1} is slightly better. Nevertheless, we present this algorithm because the ideas involved may be helpful for designing FFTs on other inverse semigroups.

\subsection{Another FFT for changing from the $\{s\}$ basis to a Fourier basis}
\label{The FFT on R_n}

Let $f \in \C R_n$ be an arbitrary element, given with respect to the $\{ s \}$ basis:
\[
f = \sum_{s\in R_n}f(s)s.
\]

As with the symmetric group (Section \ref{SnFFT}), let $t_j$ be the transposition \mbox{$(j-1,j)$}. We use the following sets of ``coset representatives" (we call them ``coset representatives" because they play the same role in this FFT as coset representatives do in group FFTs):
\begin{displaymath}
\{T_i:1\leq i \leq n, T_i = t_{i+1}t_{i+2} \ldots t_n \} \textup{ (where } T_n=Id) \textup{ and}
\end{displaymath}
\begin{displaymath}
\{T^i:1\leq i \leq n-1, T^i = t_n t_{n-1} \ldots t_{i+1} \}.
\end{displaymath}
We use the semigroup chain $R_n > R_{n-1} > \ldots > R_1$, where
\[
R_k = \{ \sigma \in R_n : \sigma(j)=j\textup{ if }j > k \}. 
\]
Halverson \cite{Halverson} has found a complete set of inequivalent, irreducible matrix representations for $R_n$ adapted to this chain. The description may be found in Section \ref{RnRepresentations}. Call this set ${\mathcal{Y}}_n$. For each $\rho \in {\Y}_n$, we must compute
\[
\hat f(\rho) = \sum_{s\in R_n}f(s)\rho(s).
\]

A subsemigroup does not necessarily partition its parent semigroup into equally sized cosets, so we cannot directly factor through a subsemigroup as in (\ref{groupfactorization}). Instead, we use an approach for $R_n$ that is based on the recursive formula (\ref{SizeRnRecursive}). With this, we have the following factorization theorem.

\begin{thm}[Factorization theorem for $R_n$]
\label{Factthm}
For any representation $\rho$ of $R_n$, if $n \geq 3$, we have the following factorization.
\begin{align}
\label{factorization}
\hat{f}(\rho) &= \sum_{i=1}^{n}\rho(T_i) \sum_{s \in R_{n-1}}f_{T_i}(s)\rho(s) +
\rho([n])\sum_{s \in R_{n-1}}f_{[n]}(s)\rho(s)
\\
&+ \sum_{i=1}^{n-1}\left[ \sum_{s\in R_{n-1}}f^{T^i}(s)\rho(s) \right] \rho(T^i), \notag{}
\end{align}
where $[n]$ is the link $(1)(2) \ldots (n-1) [n]$, $f_{A}(s) = f(As)$, and
\begin{displaymath}
f^{T^i}(s) = \begin{cases}0 & \textup{ if } n-1 \in \dom (s)
 \cr f(sT^i) & \textup{ otherwise.}\end{cases}
\end{displaymath}
\end{thm}

\begin{proof} See Section \ref{FactorizationProof}.\end{proof}

As in the group case, we will use the above breakdown to compute $\hat{f}(\rho)$ recursively. If we knew $\widehat{f_{[n]}}(\gamma)$ for all $\gamma \in {\mathcal{Y}}_{n-1}$, $\widehat{f_{T_i}}(\gamma)$ for all $\gamma \in {\mathcal{Y}}_{n-1}$ and $1\leq i \leq n$, and $\widehat{f^{T^i}}(\gamma)$ for all $\gamma \in {\mathcal{Y}}_{n-1}$ and $1\leq i \leq n-1$, using (\ref{factorization}) we could reassemble them for free, since we are using chain-adapted representations, based on how $\rho$ splits when restricted to $R_{n-1}$ (see Theorem \ref{RnBranchingThm}) to calculate $\hat{f}(\rho)$ for any $\rho \in {\mathcal{Y}}_n$. Therefore, we have
\begin{lema} \textup{For $n\geq3$, }
\[
{\mathcal{T}}_{{\mathcal Y}_n}^{s}(R_n) \leq 2n{\mathcal{T}}_{{\mathcal Y}_{n-1}}^{s}(R_{n-1}) + M_{R_n}, 
\]
where $M_{R_n}$ is the total number of operations required to compute the  sum \LP\textup{\ref{factorization}}\RP\textup{ }for all $\rho \in {\mathcal{Y}}_n$,  given knowledge of the $\widehat{f_{[n]}}$, all the $\widehat{f_{T_i}}$, and all the $\widehat{f^{T^i}}$ on $R_{n-1}$. 
\end{lema}

We now analyze $M_{R_n}$ to obtain:
\begin{thm}\textup{For $n\geq3$, }
\begin{equation}
\label{RnComplexity2}
{\mathcal T}_{{\mathcal Y}_n}^{s}(R_n) \leq 2n{\mathcal{T}}^{s}_{{\mathcal Y}_{n-1}}(R_{n-1}) + 2n^2|R_n|.
\end{equation}
\end{thm}

\begin{proof}
To analyze $M_{R_n}$, let:
\begin{itemize}
\item $M_1 =$ The maximum number of operations necessary to calculate the matrix product $\rho(T_i)A_{T_i}(\rho)$ for arbitrary matrices $A_{T_i}(\rho)$, for all $\rho \in {\mathcal{Y}}_n$ and for all $T_i$ $(1 \leq i \leq n)$.
\item $M_2 =$ The maximum number of operations necessary to calculate the matrix product $A^{T^i}(\rho)\rho(T^i)$ for arbitrary matrices $A^{T^i}(\rho)$, for all $\rho \in {\mathcal{Y}}_n$ and for all $T^i$ $(1 \leq i \leq n-1)$.
\item $M_3 =$ The maximum number of operations necessary to calculate the matrix product $\rho([n])A_{[n]}(\rho)$ for arbitrary matrices $A_{[n]}(\rho)$, for all $\rho \in {\mathcal{Y}}_n$.
\item $M_4 =$ The maximum number of operations necessary to add together $2n$ $d_\rho \times d_\rho$ arbitrary matrices, for all $\rho \in {\mathcal{Y}}_n$.
\end{itemize}
Then $M_{R_n} \leq \sum_{i=1}^4 M_i$. 

\bigskip

\noindent \underline{Analysis of $M_1$}: For each $\rho \in {\mathcal{Y}}_n$ and each $T_i$, we must perform the multiplication $\rho(T_i)A_{T_i}(\rho) = \rho(t_{i+1})\rho(t_{i+2})\ldots \rho(t_n)A_{T_i}(\rho)$ for an arbitrary matrix $A_{T_i}(\rho)$. As was the case with $S_n$, for $j >2$, $t_j \in R_j$, $t_j$ commutes with $R_{j-2}$, and ${\mathcal M}(R_j,R_{j-2}) =2$. By Schur's Lemma, $\rho(t_j)$ ($j > 2$) contains at most 2 non-zero entries per row and column. For $j=2$, $t_2 \in R_2$, and the maximum dimension of an irreducible representation of $R_2$ is 2. Therefore, $\rho(t_2)$ contains at most 2 non-zero entries per row and column. Alternatively, it is obvious from the description of ${\mathcal{Y}}_n$ (see Section \ref{RnRepresentations}) that $\rho(t_j)$ ($j \geq 2$) contains at most 2 nonzero entries per row and column. Therefore, multiplying an arbitrary matrix by $\rho(t_j)$ on the left requires at most $2d_\rho^2$ operations, and so performing the multiplication $\rho(T_i)A_{T_i}(\rho)$ requires at most $2(n-i)d_\rho^2$ operations. Therefore, we have 
\[M_1 \leq \sum_{\rho \in {\mathcal{Y}}_n}\sum_{i=1}^n 2(n-i)d_\rho^2 = (n)(n-1)|R_n|,\]
where the final equality comes from (\ref{Sum of Squares of Dimensions}).

\bigskip

\noindent \underline{Analysis of $M_2$}: The only difference between $M_1$ and $M_2$ is that $M_2$ involves multiplying arbitrary matrices by $\rho(T^i)$ on the right rather than by $\rho(T_i)$ on the left, so $M_2$ is the same as $M_1$ in the complexity analysis. Thus 
\[M_2 \leq (n)(n-1)|R_n|.\]

\bigskip

\noindent \underline{Analysis of $M_3$}: Since $[n] \in R_n$, $[n]$ commutes with $R_{n-1}$, and ${\mathcal M}(R_n, R_{n-1}) = 1$, we have that $\rho([n])$ contains at most 1 non-zero entry per row. Thus 
\[M_3 \leq \sum_{\rho \in {\mathcal{Y}}_n}d_\rho^2 = |R_n|.\]

\bigskip

\noindent \underline{Analysis of $M_4$}: For a particular $\rho$, the matrix additions can be accomplished with $(2n-1)d_\rho^2$ operations. Thus
\[M_4 \leq \sum_{\rho \in {\mathcal{Y}}_n}(2n-1)d_\rho^2 = (2n-1)|R_n|.\]
Putting this all together, we obtain
\[M_{R_n} \leq (2(n)(n-1)+1+2n-1)|R_n| = 2n^2|R_n|.\]
\end{proof}

We now prove that this algorithm gives the complexity result
\begin{thm}
For $n\geq 5$, ${\mathcal{T}}^{s}_{{\mathcal Y}_n}(R_n) \leq 2^nn|R_n|$.
\end{thm}

\begin{proof}Base case: $|R_2| = 7$, so a naive implementation of the FFT on $R_2$ gives ${\mathcal{T}}_{{\mathcal Y}_2}(R_2) \leq 49$. Applying (\ref{RnComplexity2}) repeatedly, we have
\begin{align*}
{\mathcal{T}}_{{\mathcal Y}_3}(R_3) &\leq 2(3){\mathcal{T}}_{{\mathcal Y}_2}(R_2) + 2(3)^2|R_3| \leq 6(49) + 18(34) = 906, \textup{ so } \\
{\mathcal{T}}_{{\mathcal Y}_4}(R_4) &\leq 2(4){\mathcal{T}}_{{\mathcal Y}_3}(R_3) + 2(4)^2|R_4| \leq 8(906) + 32(209) = 13936\textup{, so} \\
{\mathcal{T}}_{{\mathcal Y}_5}(R_5)  &\leq 2(5){\mathcal{T}}_{{\mathcal Y}_4}(R_4) + 2(5)^2|R_5| \leq 10(13936) + 50(1546) = 216660 \textup{, and } \\
& 216660 < 2^5(5)|R_5| = 247360.  \end{align*}
This proves the base case.

Similarly, for $n=6$, we have
\begin{align*}
{\mathcal{T}}_{{\mathcal Y}_6}(R_6) & \leq 2(6){\mathcal{T}}_{{\mathcal Y}_5}(R_5) + 2(6)^2|R_6| \leq 12(216660) + 72(13327) \\
&= 3559464,\textup{ and } 3559464 < 2^6(6)|R_6| = 5117568.\end{align*}

Now, let $n \geq 7$. Observe that, for $\alpha = 1 \textup{ to } n$, the sets  $\{\sigma \in R_n : \sigma(\alpha)=n\}$ are disjoint, and each are of size $|R_{n-1}|$; thus $n|R_{n-1}| \leq |R_n|$. Therefore, we have
\begin{align*}
{\mathcal{T}}_{{\mathcal Y}_n}(R_n) &\leq 2n{\mathcal{T}}_{{\mathcal Y}_{n-1}}(R_{n-1}) + 2n^2|R_n| \\ &\leq
2n(2^{n-1}(n-1)|R_{n-1}|) + 2n^2|R_n| \\ &\leq 
2^n(n-1)|R_n| + 2n^2|R_n| \\ &=
2^nn|R_n| + (2n^2 - 2^n)|R_n| \\ &\leq
2^nn|R_n| \textup{ (since } 2n^2 \leq 2^n \textup{ for }n \geq 7).
\end{align*}
\end{proof}


\subsection{Chain-adapted matrix representations for $R_n$}
\label{RnRepresentations}
In this section, we give a description of a complete set of irreducible, inequivalent, chain-adapted matrix representations for $R_n$ relative to the chain $R_n > R_{n-1} > \ldots > R_1$. The results herein are a special case of the results in \cite{Halverson}.

\begin{defa}[partition] A {\em partition} $\lambda$ of a nonnegative integer $k$ \LP written $\lambda \vdash k$\RP\textup{ }is a weakly decreasing sequence of nonnegative integers whose sum is $k$. We consider two partitions to be equal if and only if they only differ by the number of $0$'s they contain, and we identify a partition $\lambda$ with its Young diagram. \end{defa}

For example, $\lambda = (5,5,3,1)$ is a partition of $14$, and  
\[\lambda = (5,5,3,1) = (5,5,3,1,0) = \]
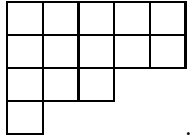
\[
{\begin{Young}
     \;  &  \; &  \; & \;  &\; \cr
\;  &  \; &  \; & \;  &\; \cr
    \;  &  \;  & \;  \cr
    \; \cr
     \end{Young}}.\]

It is well-known that a complete set of inequivalent, irreducible representations for $R_n$ is indexed by the set of all partitions of the integers $\{0, 1, \ldots n\}$ (see, for example, \cite{Grood} or \cite{Solomon}). Therefore, for integers $n \geq 0$, let
\begin{displaymath}
\Lambda_n = \{\lambda \vdash k : 0 \leq k \leq n\}.
\end{displaymath}

\begin{defa}[n-tableau, n-standard tableau]For $\lambda \in \Lambda_n$, define $L$ to be an {\em $n$-tableau of shape $\lambda$} if it is a filling of the diagram for $\lambda$ with numbers from $\{1, 2, \ldots n\}$ such that each number in $L$ appears at most once. $L$ is an {\em $n$-standard tableau} if, furthermore, the entries in each column of $L$ increase from top to bottom and the entries in each row of $L$ increase from left to right.
\end{defa}

Fix $\lambda$. Let ${T}_n^\lambda$ denote the set of $n$-standard tableaux of shape $\lambda$. The symmetric group acts on tableaux by permuting their entries. If L is an $n$-tableau, then $(i-1,i)L$ is the tableau that is obtained from $L$ by replacing $i-1$ (if $i-1 \in L$) with $i$, and by replacing $i$ (if $i \in L$) with $i-1$. Note that $L \in {T}_n^\lambda$ need not imply $(i-1, i)L \in {T}_n^\lambda$.

Let $\{v_L : L \in {T}_n^\lambda\}$ be a set of independent vectors. We form
\begin{displaymath}
V^\lambda = \mathbb C\textup{-span}\{v_L : L \in T_n^\lambda\}.
\end{displaymath}
As such, the symbols $v_L$, for $L\in T_n^\lambda$, are a basis for the vector space $V^\lambda$. Halverson defines an action of $R_n$ on $V^\lambda$ in such a way that (extending by linearity) $V^\lambda$ is an irreducible $\mathbb CR_n$-module and such that, as $\lambda$ ranges over $\Lambda_n$, the $V^\lambda$ constitute a complete set of inequivalent, irreducible representation modules for $R_n$. We first describe this action, and we then describe an ordering of the bases for the $V^\lambda$ so that the resulting matrix representations are chain-adapted to $R_n > R_{n-1} > \ldots > R_1$.


\begin{defa}[content] If $b$ is a box of $\lambda$ in position $(i,j)$, then the {\em content} of $b$ is defined to be
\begin{displaymath}
ct(b) = j-i.
\end{displaymath}
\end{defa}

Let $L \in T^\lambda_n$. If $i-1, i \in L$, then let $L(i-1)$ and $L(i)$ denote the box in $L$ containing $i-1$ and $i$, respectively.

To define the action of $R_n$ on $V^\lambda$, it is sufficient to define the action of a set of generators of $R_n$ on $V^\lambda$.

\begin{defa}[action of $R_n$ on $V^\lambda$] Define the action of the transpositions $t_i = (i-1, i)$, for $2\leq i \leq n$, as follows:
\begin{displaymath}
t_i v_L = \begin{cases} 
\frac{1}{ct(L(i))-ct(L(i-1))}v_L + (1+\frac{1}{ct(L(i))-ct(L(i-1))})v_{L'} & \textup{if }i-1,i \in L \cr
v_{t_iL} & \textup{if exactly one of} 
\cr & i-1,i \in L \cr
v_L & \textup{if } i-1,i \notin L
\end{cases}
\end{displaymath}
where
\begin{displaymath}
v_{L'} = \begin{cases}
v_{t_iL} & \textup{if } t_iL \textup{ is } n\textup{-standard} \cr
0 & \textup{otherwise.} \end{cases}
\end{displaymath}

Define the action of the link $(1)(2)\ldots(n-1)[n] = [n]$ on $V^\lambda$ by
\begin{displaymath}
[n]v_L = \begin{cases} v_L & \textup{ if } n \notin L \cr
0 & \textup{ if } n \in L.
\end{cases}
\end{displaymath}
\end{defa}

\bigskip \noindent {\bf Remark}: If $\lambda=(0)$, then $V^\lambda$ is 1-dimensional, and the action of $R_n$ on $V^\lambda$ is the trivial action given by $xv=v$ for all $x\in R_n$ and all $v \in V^\lambda$.

\begin{thm} As $\lambda$ varies over all partitions of all nonnegative integers less than or equal to $n$, the $V^\lambda$ constitute a complete set of irreducible, pairwise non-isomorphic representation modules for $R_n$ \cite{Halverson}.
\end{thm}

\begin{defa}[corner of a partition] A corner is a box $c$ of $\lambda$ for which $\lambda$ contains no box to the right or below $c$. In other words, the corners are the possible positions of $n$ in an $n$-standard tableau of shape $\lambda$.
\end{defa}

We now record the branching theorem for $R_n$.
\begin{thm}[Branching theorem]
As a $\mathbb CR_{n-1}$ module, 
\begin{displaymath}
V^\lambda \cong \oplus_{\mu \in \lambda^{-,=}}V^\mu,
\end{displaymath}
where $\lambda^{-,=}$ is the set of all partitions $\mu \in \Lambda_{n-1}$ such that either $\mu = \lambda$ \LP if $\lambda \nvdash n$\RP\textup{ }or $\mu$ is obtained by removing a corner from $\lambda$ \cite{Halverson}.
\end{thm}

Now, for purposes of chain-adaptation, we order the basis $\{v_L\}$ for $V^\lambda$ using the following generalized last-letter ordering.

We begin by partitioning the $v_L$ into subsets based on the corners $c_1, \ldots, c_l$ of $\lambda$. Number the corners from top to bottom. Now, form the sets
\begin{displaymath}
V^\lambda(0)=\{v_L : L \in T^\lambda_n \textup{ and n} \notin L\},
\end{displaymath}
\begin{displaymath}
V^\lambda(i)=\{v_L : L \in T^\lambda_n \textup{ and n} \in c_i \textup{ of }L\} , \quad 1\leq i\leq l,
\end{displaymath}
and declare elements of $V^\lambda(j)$ to be earlier in the ordering than elements of $V^\lambda(k)$ whenever $j < k$. To order the subset $V^\lambda(k)$, delete the corner $c_k$ (do nothing if $k=0$) and repeat the same ordering process (starting by identifying the corners of the resulting partition and partitioning the $v_L$ into subsets based on those corners) with $n-1$ in place of $n$, etc.

As an example, consider $\lambda = (2,1,1)$, which has two corners, and $R_5$. Our ordered basis for the 15-dimensional $V^\lambda$ is
\begin{small}
\begin{align*}
v_{\begin{Young}1&4 \cr 2 \cr 3 \cr \end{Young}} &<
v_{\begin{Young}1&3 \cr 2 \cr 4 \cr \end{Young}} <
v_{\begin{Young}1&2 \cr 3 \cr 4 \cr \end{Young}} <   
v_{\begin{Young}1&5 \cr 2 \cr 3 \cr \end{Young}} <
v_{\begin{Young}1&5 \cr 2 \cr 4 \cr \end{Young}} <
v_{\begin{Young}1&5 \cr 3 \cr 4 \cr \end{Young}} <
v_{\begin{Young}2&5 \cr 3 \cr 4 \cr \end{Young}} <   
v_{\begin{Young}1&3 \cr 2 \cr 5 \cr \end{Young}} 
\\
&< v_{\begin{Young}1&2 \cr 3 \cr 5 \cr \end{Young}} <
v_{\begin{Young}1&4 \cr 2 \cr 5 \cr \end{Young}} <
v_{\begin{Young}1&4 \cr 3 \cr 5 \cr \end{Young}} <
v_{\begin{Young}2&4 \cr 3 \cr 5 \cr \end{Young}} <
v_{\begin{Young}1&2 \cr 4 \cr 5 \cr \end{Young}} <
v_{\begin{Young}1&3 \cr 4 \cr 5 \cr \end{Young}} <
v_{\begin{Young}2&3 \cr 4 \cr 5 \cr \end{Young}}.
\end{align*}
\end{small}

\bigskip \noindent {\bf Remark}: If $\lambda \vdash n$, then the generalized last-letter ordering scheme given above reduces to the usual last-letter ordering scheme used for Young's orthogonal and seminormal representations of the symmetric group.

It is now easy to see, under this ordering of the bases for the $V^\lambda$, that the matrix representations described in this section are chain-adapted to the chain \mbox{$R_n > R_{n-1} > \ldots > R_1$}. We know, by the branching theorem for $R_n$, how $V^\lambda$ decomposes as a $R_{n-1}$ module (it decomposes into the $\C$-span of the $V^\lambda(k)$ for \mbox{$0 \leq k \leq l$}), and it is obvious by the action of $R_{n-1}$ on $V^\lambda$ that 
\[
R_{n-1} \C\textup{-span}(V^\lambda(k)) \subseteq \C\textup{-span}(V^\lambda(k))
\]
for all $k$. The same argument is used to induct down the chain $R_n > R_{n-1} > \ldots > R_1$, now, and is trivial because we ordered our basis for $V^\lambda$ inductively according to the same rule. We now restate the branching theorem for $R_n$ under our ordering of the bases for the $V^\lambda$.

\begin{thm}[Branching theorem]
\label{RnBranchingThm}
Let $\rho^\lambda$ be the matrix representation associated to $V^\lambda$ with respect to the basis $\{v_L\}$, with the basis ordered according to the generalized last-letter ordering. Then 
\begin{displaymath}
\rho^\lambda|_{R_{n-1}} = \oplus_{\mu \in \lambda^{-,=}}\rho^\mu,
\end{displaymath}
where $\lambda^{-,=}$ is the set of all partitions $\mu \in \Lambda_{n-1}$ such that either $\mu = \lambda$ or $\mu$ is obtained by removing a corner from $\lambda$. The first $\mu \in \lambda^{-,=}$ is the one that removes no corners \LP if $\lambda \nvdash n$\RP, the next $\mu$ is the one that removes the highest corner, the next $\mu$ is the one that removes the second highest corner, etc. 
\end{thm}


\subsection{Proof of Theorem \ref{Factthm}}
\label{FactorizationProof}

To provide motivation for the proof of Theorem \ref{Factthm}, we begin by proving the recursive formula
\begin{thma}[Theorem \ref{SizeRnRecursiveThm}]For $n \geq 3$, $|R_n| = 2n|R_{n-1}| - (n-1)^2|R_{n-2}|$. \end{thma}

\begin{proof}
Viewing the elements of $R_n$ as rook matrices, $R_n$ consists of those elements having all 0's in column $n$ and row $n$ (of which there are $|R_{n-1}|$), together with, for each $\alpha \in \{1,\ldots, n\}$, those having a 1 in position $(\alpha, n)$ (of which there are $n|R_{n-1}|$ total), together with, for each $\alpha \in \{1,\ldots, n-1\}$, those having a 1 in position $(n,\alpha)$ (of which there are $(n-1)|R_{n-1}|$ total). Counting the number of elements of $R_n$ in this way overcounts. For each $(\alpha,\beta)$ with $1 \leq \alpha , \beta \leq n-1$, every element with 1's in positions $(\alpha,n)$ and $(n,\beta)$ (of which there are $(n-1)^2|R_{n-2}|$ total) gets counted twice.\end{proof}

Now, let $f\in \C R_n$ be given with respect to the $\{s\}$ basis. We prove Theorem \ref{Factthm}.

\begin{thma}[Theorem \ref{Factthm}]For any representation $\rho$ on $R_n$, if $n \geq 3$, we have the following factorization.
\begin{align*}
\hat{f}(\rho) & = \sum_{i=1}^{n}\rho(T_i) \sum_{s \in R_{n-1}}f_{T_i}(s)\rho(s) +
\rho([n])\sum_{s \in R_{n-1}}f_{[n]}(s)\rho(s) \\
&+ \sum_{i=1}^{n-1}\left[ \sum_{s\in R_{n-1}}f^{T^i}(s)\rho(s) \right] \rho(T^i),
\end{align*}
where $[n]$ is the link $(1)(2) \ldots (n-1) [n]$, $f_{A}(s) = f(As)$, and
\begin{displaymath}
f^{T^i}(s) = \begin{cases}0 & \textup{ if } n-1 \in \dom (s)
 \\ f(sT^i) & \textup{ otherwise.}\end{cases}
\end{displaymath}
\end{thma}

\begin{proof}
Let $n \geq 3$. We have 3 types of elements $s$ of $R_n$.
\begin{itemize}
\item Type 1: Those for which $s(n)=i$ for some $1 \leq i \leq n$.
\item Type 2: Those for which both $s(i)=n$ for some $1 \leq i \leq n-1$ and \mbox{$n \notin$ dom($s$)}.
\item Type 3: Those for which both $s(i) \neq n$ for all $1 \leq i \leq n$ and $n \notin$ dom($s$).
\end{itemize}
By the argument given in the proof of Theorem \ref{SizeRnRecursiveThm}, this counts all elements of $R_n$ precisely once. 

Now, let $1 \leq i \leq n$. View the ``coset representative" $T_i$ as a permutation matrix, and view the elements $s \in R_n$ as rook matrices. Multiplying any matrix $X$ on the left by $T_i$ simply moves row $j$ of $X$ to row $j+1$ (for all $j$ such that $i \leq j \leq n-1$) and moves row $n$ of $X$ to row $i$. Thus, as $s$ varies over $R_{n-1}$, $T_is$ varies bijectively over $\{s \in R_n : s(n) = i\}$. 
Therefore, we have
\begin{align*}
\sum_{s \in R_n \textup{ of Type 1}}f(s)\rho(s) &=
\sum_{i=1}^n \sum_{s \in R_{n-1}} f(T_is)\rho(T_is) \\ &=
\sum_{i=1}^n\sum_{s \in R_{n-1}}f_{T_i}(s)\rho(T_i)\rho(s) \\ &=
\sum_{i=1}^n \rho(T_i) \sum_{s \in R_{n-1}}f_{T_i}(s)\rho(s),
\end{align*}
where $f_{T_i}(s) = f(T_is)$.

Similarly, multiplying any matrix $X$ on the right by $T^i$ moves column $j$ of $X$ to column $j+1$ $(i \leq j \leq n-1)$ and moves column $n$ of $X$ to column $i$. Thus, as $s$ varies over $R_{n-1}$, $sT^i$ varies bijectively over $\{s \in R_n : s(i) = n\}$. So
\begin{align*}
\sum_{s \in R_n:s(i)=n} f(s)\rho(s) & = 
\sum_{s \in R_{n-1}} f(sT^i)\rho(sT^i) \\ &=
\sum_{s \in R_{n-1}} f(sT^i)\rho(s)\rho(T^i).
\end{align*}
To ensure that we only count the elements of Type 2, we restrict our attention to $1 \leq i \leq n-1$, and we set the function values of the elements of Type 1 appearing in the above sum to 0:
\begin{align*}
\sum_{s \in R_n \textup{ of Type 2}}f(s)\rho(s) &=
\sum_{i=1}^{n-1} \sum_{s \in R_{n-1}} f^{T^i}(s)\rho(sT^i) \\ &= 
\sum_{i=1}^{n-1} \left[ \sum_{s \in R_{n-1}} f^{T^i}(s)\rho(s)\right] \rho(T^i),
\end{align*}
where \[f^{T^i}(s) = \begin{cases}0 & \textup{ if } n-1 \in \dom (s) \textup{ (i.e. } n \in \dom (sT^i))
\cr f(sT^i) & \textup{ otherwise.}\end{cases}\]

Finally, 
\begin{align*}
\sum_{s \in R_n \textup{ of Type 3}}f(s)\rho(s) &=
\sum_{s \in R_{n-1}}f([n]s)\rho([n]s) \\ &= 
\rho([n])\sum_{s \in R_{n-1}}f_{[n]}(s)\rho(s).
\end{align*}
Putting this all together, then, we find that for any representation $\rho$ of $R_n$ and any $n \geq 3$,
\begin{align*}
\hat{f}(\rho) &= 
\sum_{s \in R_n \textup{ of Type 1}} f(s)\rho(s) + \sum_{s \in R_n \textup{ of Type 2}} f(s)\rho(s) + \sum_{s \in R_n \textup{ of Type 3}} f(s)\rho(s) \\ &=
\sum_{i=1}^{n}\rho(T_i) \sum_{s \in R_{n-1}}f_{T_i}(s)\rho(s) +
\rho([n])\sum_{s \in R_{n-1}}f_{[n]}(s)\rho(s) \\
&+ \sum_{i=1}^{n-1}\left[ \sum_{s\in R_{n-1}}f^{T^i}(s)\rho(s) \right] \rho(T^i).
\end{align*}
\end{proof}

\begin{center}
{\bf Concluding remarks}
\end{center}

The extension of FFTs to semigroups creates a new collection of interesting challenges. We remark that many of the ideas in this paper (such as those from \cite{Steinberg2} and several of the results in this paper that follow from them) can be extended to general inverse semigroups. As with groups, FFT algorithms for inverse semigroups will vary from semigroup to semigroup, but a number of the underlying ideas are the same for any inverse semigroup FFT. These ideas can be cast in a general framework, which we intend to help guide the development of future FFTs, and this general framework is the subject of a paper currently in preparation.

Also, as mentioned previously, spectral analysis for the rook monoid involves projecting a function onto the isotypic subspaces of $\C R_n$, which can be accomplished by means of an FFT, and examining the resulting projections. A variety of interesting issues arise in this analysis, and we have worked out a detailed example, consisting of partially ranked voting data on $R_5$, to explain them. These results are the subject of a paper currently in preparation.

\nocite{DanCharacterProjection}
\nocite{Dan1990}
\nocite{DanSepVars}
\nocite{DanSchurs}
\nocite{CooleyTukey}
\nocite{Lawson}
\nocite{James}
\nocite{Steinberg2}
\nocite{Steinberg}
\nocite{Persi}
\nocite{Clausen}
\nocite{Solomon}
\nocite{Grood}
\nocite{Halverson}
\nocite{Holmes}
\nocite{Asymptotics}
\nocite{DanDiameters}
\nocite{CliffandPres}
\nocite{CliffandPres2}
\nocite{Munn1}
\nocite{Munn2}
\nocite{Munn3}
\nocite{Rhodes}
\nocite{Serre}
\nocite{DennisFarb}
\nocite{DanAMS}



\begin{thebibliography}{10}

\bibitem{Baum}
U.~Baum, \emph{Existence and efficient construction of fast {F}ourier
  transforms for supersolvable groups}, Comput. Complexity \textbf{1} (1991),
  235--256.

\bibitem{Clausen}
M.~Clausen and U.~Baum, \emph{Fast {F}ourier transforms for symmetric groups:
  Theory and implementation}, Math. Comput. \textbf{61}, no.~204.

\bibitem{Clifford}
A.~H. Clifford, \emph{{M}atrix representations of completely simple
  semigroups}, Amer. J. Math. \textbf{64} (1942), no.~1/4, 327--342.

\bibitem{CliffandPres}
A.~H. Clifford and G.~B. Preston, \emph{The {A}lgebraic {T}heory of
  {S}emigroups}, vol.~1, Mathematical {S}urveys {N}o. 7, {AMS}, {P}rovidence,
  {RI}, 1961.

\bibitem{CliffandPres2}
\bysame, \emph{The {A}lgebraic {T}heory of {S}emigroups}, vol.~2, Mathematical
  {S}urveys {N}o. 7, {AMS}, {P}rovidence, {RI}, 1961.

\bibitem{CooleyTukey}
J.~W. Cooley and J.~W. Tukey, \emph{An algorithm for machine calculation of
  complex {F}ourier series}, Math. Comput. \textbf{19} (1965), 297--301.

\bibitem{Curtis}
C.~Curtis and I.~Reiner, \emph{Representation {T}heory of {F}inite {G}roups and
  {A}ssociative {A}lgebras}, John Wiley and Sons, 1962.

\bibitem{Persi}
P.~Diaconis, \emph{A generalization of spectral analysis with application to
  ranked data}, Ann. Statist. \textbf{17} (Sept. 1989), no.~3.

\bibitem{Dan1990}
P.~Diaconis and D.~Rockmore, \emph{Efficient computation of the {F}ourier
  transform on finite groups}, J. Amer. Math. Soc. \textbf{3} (April, 1990),
  no.~2, 297--332.

\bibitem{DanCharacterProjection}
\bysame, \emph{Efficient computation of isotypic projections for the symmetric
  group}, DIMACS Series in Discrete Mathematics and Theoretical Computer
  Science \textbf{00} (March 1, 1992), no.~0000.

\bibitem{DennisFarb}
B.~Farb and R.~K. Dennis, \emph{Noncommutative algebra}, Graduate Texts in
  Mathematics, vol. 144, Springer-Verlag, New York-Heidelberg, 1993.

\bibitem{Grood}
C.~Grood, \emph{A {S}pecht module analog for the rook monoid}, Electron. J.
  Combin. \textbf{9} (2002).

\bibitem{Halverson}
T.~Halverson, \emph{Representations of the q-rook monoid}, J. Algebra
  \textbf{273} (2004), 227--251.

\bibitem{Holmes}
R.~B. Holmes, \emph{Mathematical foundation of signal processing {II}. {T}he
  role of group theory}, Massachusetts Institute of Technology Technical Report
  781 (October 13, 1987).

\bibitem{James}
G.~D. James, \emph{The representation theory of the symmetric groups}, Lect.
  Notes Math., Springer-Verlag, Berlin \textbf{682} (1978).

\bibitem{Asymptotics}
S.~Janson and V.~Mazorchuk, \emph{Some remarks on the combinatorics of
  {$IS_n$}}, Semigroup Forum \textbf{70} (June 2, 2005), no.~3, 391--405.

\bibitem{Lawson}
M.~V. Lawson, \emph{Inverse semigroups: The theory of partial symmetries},
  World Scientific, Singapore, 1998.

\bibitem{Maslen}
D.~K. Maslen, \emph{The efficient computation of {F}ourier transforms on the
  symmetric group}, Math. Comput. \textbf{67} (1998), no.~223, 1121--1147.

\bibitem{DanSchurs}
D.~K. Maslen and D.~N. Rockmore, \emph{Generalized {FFT}s - a survey of some
  recent results}, Proceedings of the DIMACS Workshop on Groups and Computation
  (1997).

\bibitem{DanSepVars}
\bysame, \emph{Separation of variables and the computation of {F}ourier
  transforms on finite groups, {I}}, J. Amer. Math. Soc. \textbf{10} (January
  1997), no.~1, 169--214.

\bibitem{DanAMS}
\bysame, \emph{The {C}ooley-{T}ukey {FFT} and group theory}, Notices of the AMS
  \textbf{48} (November 2001), no.~10, 1151--1161.

\bibitem{Munn1}
W.~D. Munn, \emph{On semigroup algebras}, Proc. Cambridge Philos. Soc.
  \textbf{51} (1955), 1--15.

\bibitem{Munn3}
\bysame, \emph{The characters of the symmetric inverse semigroup}, Proc.
  Cambridge Philos. Soc. \textbf{53} (1957), 13--18.

\bibitem{Munn2}
\bysame, \emph{Matrix representations of semigroups}, Proc. Cambridge Philos.
  Soc. \textbf{53} (1957), 5--12.

\bibitem{Rhodes}
J.~Rhodes and Y.~Zalcstein, \emph{Monoids and semigroups with applications},
  ch.~title: {E}lementary representation and character theory of finite
  semigroups and its application, pp.~334--367, World Sci. Publishing, River
  Edge, NJ, 1991.

\bibitem{Dan2004}
D.~Rockmore, \emph{{NATO} science series, computational noncommutative algebra
  and applications}, vol. 136, ch.~title: {R}ecent Progress and Applications in
  Group {FFT}s, pp.~227--254, Springer Netherlands, 2004.

\bibitem{DanDiameters}
D.~N. Rockmore and D.~K. Maslen, \emph{Adapted {D}iameters and {FFT}s on
  {G}roups}, Proc. 6 th ACM-SIAM SODA, 253--262.

\bibitem{Serre}
J.~P. Serre, \emph{Linear representations of finite groups}, Graduate Texts in
  Mathematics, vol.~42, Springer-Verlag, New York-Heidelberg, 1977.

\bibitem{Solomon}
L.~Solomon, \emph{Representations of the rook monoid}, J. Algebra \textbf{256}
  (2002), 309--342.

\bibitem{Stanley}
R.~Stanley, \emph{Enumerative combinatorics. vol. 1}, Cambridge Studies in
  Advanced Mathematics, vol.~49, Cambridge University Press, 1997.

\bibitem{Steinberg2}
B.~Steinberg, \emph{Mobius functions and semigroup representation theory {II}:
  Character formulas and multiplicities}, Preprint.

\bibitem{Steinberg}
\bysame, \emph{Mobius functions and semigroup representation theory}, J. Comb.
  Theor. A. \textbf{113} (2006), 866--881.

\bibitem{Yates}
F.~Yates, \emph{The design and analysis of factorial experiments}, Imp. Bur.
  Soil Sci. Tech. Comm. \textbf{35} (1937).

\end{thebibliography}

\providecommand{\bysame}{\leavevmode\hbox to3em{\hrulefill}\thinspace}
\providecommand{\MR}{\relax\ifhmode\unskip\space\fi MR }
\providecommand{\MRhref}[2]{%
  \href{http://www.ams.org/mathscinet-getitem?mr=#1}{#2}
}
\providecommand{\href}[2]{#2}


\end{document}